\documentclass{amsart}
\usepackage{amssymb}

\usepackage{amsmath}
\usepackage{amscd}
\usepackage{thmdefs}

\input{tcilatex}
\theoremstyle{definition}
\theoremstyle{remark}
\numberwithin{equation}{section}

\input tcilatex

\begin{document}
\title[Random Variables in a Graph $W^{*}$-Probability Space]{Random Variables in a Graph $W^{*}$-Probability Space}
\author{Ilwoo Cho}
\address{Univ. of Iowa, Dep. of Math, Iowa City, IA, U. S. A}
\email{ilcho@math.uiowa.edu}
\date{}
\dedicatory{}
\thanks{}
\maketitle

\begin{abstract}
In [15], we constructed a $W^{*}$-probability space $\left(
W^{*}(G),E\right) $ with amalgamation over a von Neumann algebra $D_{G},$
where $W^{*}(G)$ is a graph $W^{*}$-algebra induced by the countable
directed graph $G$. In [15], we computed the $D_{G}$-valued moments and
cumulants of arbitrary random variables in $\left( W^{*}(G),E\right) $ and
we could characterize the $D_{G}$-freeness of generators of $W^{*}(G),$ by
the so-called diagram-distinctness on $G.$ In this paper, we will observe
some special $D_{G}$-valued random variables in $\left( W^{*}(G),E\right) ,$
for instance, $D_{G}$-semicircular elements, $D_{G}$-even elements, $D_{G}$%
-valued R-diagonal elements and the generating operator of $W^{*}(G).$ In
particular, we can get that (i) if $l$ is a loop in the graph $G,$ then the
random variable $L_{l}+L_{l}^{*}$ is $D_{G}$-semicircular, (ii) if $w$ is a
finite path, then the random variable $L_{w}+L_{w}^{*}$ is $D_{G}$-even,
(iii) if $w$ is a finite path, then the random variables $L_{w}$ and $%
L_{w}^{*}$ are $D_{G}$-valued R-diagonal.
\end{abstract}

\strut

In this paper, we construct the graph $W^{*}$-probability spaces. The graph $%
W^{*}$-probability theory is one of the good example of Speicher's
combinatorial free probability theory with amalgamation (See [16]). In this
paper, we will observe how to compute the moment and cumulant of an
arbitrary random variables in the graph $W^{*}$-probability space and the
freeness on it with respect to the given conditional expectation. Also, we
consider certain special random variables of the graph $W^{*}$-probability
space, for example, semicircular elements, even elements and R-diagonal
elements. This shows that the graph $W^{*}$-probability spaces contain the
rich free probability objects.

\bigskip

In [10], Kribs and Power defined the free semigroupoid algebras and obtained
some properties of them. Our work is highly motivated by [10]. Roughly
speaking, graph $W^{*}$-algebras are $W^{*}$-topology closed version of free
semigroupoid algebras. Throughout this paper, let $G$ be a countable
directed graph and let $\mathbb{F}^{+}(G)$ be the free semigroupoid of $G,$
in the sense of Kribs and Power. i.e., it is a collection of all vertices of
the graph $G$ as units and all admissible finite paths, under the
admissibility. As a set, the free semigroupoid $\mathbb{F}^{+}(G)$ can be
decomposed by

\strut

\begin{center}
$\mathbb{F}^{+}(G)=V(G)\cup FP(G),$
\end{center}

\strut

where $V(G)$ is the vertex set of the graph $G$ and $FP(G)$ is the set of
all admissible finite paths. Trivially the edge set $E(G)$ of the graph $G$
is properly contained in $FP(G),$ since all edges of the graph can be
regarded as finite paths with their length $1.$ We define a graph $W^{*}$%
-algebra of $G$ by

\strut

\begin{center}
$W^{*}(G)\overset{def}{=}\overline{%
\mathbb{C}[\{L_{w},L_{w}^{*}:w\in
\mathbb{F}^{+}(G)\}]}^{w},$
\end{center}

\strut

where $L_{w}$ and $L_{w}^{\ast }$ are creation operators and annihilation
operators on the generalized Fock space $H_{G}=l^{2}\left( \mathbb{F}%
^{+}(G)\right) $ induced by the given graph $G,$ respectively. Notice that
the creation operators induced by vertices are projections and the creation
operators induced by finite paths are partial isometries. We can define the $%
W^{\ast }$-subalgebra $D_{G}$ of $W^{\ast }(G),$ which is called the
diagonal subalgebra by

\strut

\begin{center}
$D_{G}\overset{def}{=}\overline{\mathbb{C}[\{L_{v}:v\in V(G)\}]}^{w}.$
\end{center}

\strut

Then each element $a$ in the graph $W^{*}$-algebra $W^{*}(G)$ is expressed by

\strut

\begin{center}
$a=\underset{w\in \mathbb{F}^{+}(G:a),\,u_{w}\in \{1,*\}}{\sum }%
p_{w}L_{w}^{u_{w}},$ \ for $p_{w}\in \mathbb{C},$
\end{center}

\strut

where $\mathbb{F}^{+}(G:a)$ is a support of the element $a$, as a subset of
the free semigroupoid $\mathbb{F}^{+}(G).$ The above expression of the
random variable $a$ is said to be the Fourier expansion of $a.$ Since $%
\mathbb{F}^{+}(G)$ is decomposed by the disjoint subsets $V(G)$ and $FP(G),$
the support $\mathbb{F}^{+}(G:a)$ of $a$ is also decomposed by the following
disjoint subsets,

\strut

\begin{center}
$V(G:a)=\mathbb{F}^{+}(G:a)\cap V(G)$
\end{center}

and

\begin{center}
$FP(G:a)=\mathbb{F}^{+}(G:a)\cap FP(G).$
\end{center}

\strut

Thus the operator $a$ can be re-expressed by

\strut

\begin{center}
$a=\underset{v\in V(G:a)}{\sum }p_{v}L_{v}+\underset{w\in FP(G:a),\,u_{w}\in
\{1,*\}}{\sum }p_{w}L_{w}^{u_{w}}.$
\end{center}

\strut

Notice that if $V(G:a)\neq \emptyset ,$ then $\underset{v\in V(G:a)}{\sum }%
p_{v}L_{v}$ is contained in the diagonal subalgebra $D_{G}.$ Thus we have
the canonical conditional expectation $E:W^{*}(G)\rightarrow D_{G},$ defined
by

\strut

\begin{center}
$E\left( a\right) =\underset{v\in V(G:a)}{\sum }p_{v}L_{v},$
\end{center}

\strut

for all $a=\underset{w\in \mathbb{F}^{+}(G:a),\,u_{w}\in \{1,*\}}{\sum }%
p_{w}L_{w}^{u_{w}}$ \ in $W^{*}(G).$ Then the algebraic pair $\left(
W^{*}(G),E\right) $ is a $W^{*}$-probability space with amalgamation over $%
D_{G}$ (See [16]). It is easy to check that the conditional expectation $E$
is faithful in the sense that if $E(a^{*}a)=0_{D_{G}},$ for $a\in W^{*}(G),$
then $a=0_{D_{G}}.$

\strut

For the fixed operator $a\in W^{*}(G),$ the support $\mathbb{F}^{+}(G:a)$ of
the operator $a$ is again decomposed by

\strut

\begin{center}
$\mathbb{F}^{+}(G:a)=V(G:a)\cup FP_{*}(G:a)\cup FP_{*}^{c}(G:a),$
\end{center}

\strut

with the decomposition of $FP(G:a),$

$\strut $

\begin{center}
$FP(G:a)=FP_{*}(G:a)\cup FP_{*}^{c}(G:a),$
\end{center}

where

\strut

\begin{center}
$FP_{*}(G:a)=\{w\in FP(G:a):$both $L_{w}$ and $L_{w}^{*}$ are summands of $%
a\}$
\end{center}

and

\begin{center}
$FP_{*}(G:a)=FP(G:a)\,\,\setminus \,\,FP_{*}(G:a).$
\end{center}

\strut

The above new expression plays a key role to find the $D_{G}$-valued moments
of the random variable $a.$ In fact, the summands $p_{v}L_{v}$'s and $%
p_{w}L_{w}+p_{w^{t}}L_{w}^{*},$ for $v\in V(G:a)$ and $w\in FP_{*}(G:a)$ act
for the computation of $D_{G}$-valued moments of $a.$ By using the above
partition of the support of a random variable, we can compute the $D_{G}$%
-valued moments and $D_{G}$-valued cumulants of it via the lattice path
model $LP_{n}$ and the lattice path model $LP_{n}^{*}$ satisfying the $*$%
-axis-property. At a first glance, the computations of $D_{G}$-valued
moments and cumulants look so abstract (See Chapter 3) and hence it looks
useless. However, these computations, in particular the computation of $%
D_{G} $-valued cumulants, provides us how to figure out the $D_{G}$-freeness
of random variables by making us compute the mixed cumulants. As
applications, in the final chapter, we can compute the moment and cumulant
of the operator that is the sum of $N$-free semicircular elements with their
covariance $2.$ If $a$ is the operator, then the $n$-th moment of $a$ is

\strut

\begin{center}
$\left\{ 
\begin{array}{ll}
\left( 2N\right) ^{\frac{n}{2}}\cdot c_{\frac{n}{2}} & \text{if }n\text{ is
even} \\ 
0 & \text{if }n\text{ is odd,}
\end{array}
\right. $
\end{center}

\strut

and the $n$-th cumulant of $a$ is

\strut

\begin{center}
$\left\{ 
\begin{array}{lll}
2N &  & \text{if }n=2 \\ 
0 &  & \text{otherwise.}
\end{array}
\right. $
\end{center}

\strut \strut \strut

\begin{center}
Based on the $D_{G}$-cumulant computation, we can characterize the $D_{G}$%
-freeness of generators of $W^{*}(G),$ by the so-called diagram-distinctness
on the graph $G.$ i.e., the random variables $L_{w_{1}}$ and $L_{w_{2}}$ are
free over $D_{G}$ if and only if $w_{1}$ and $w_{2}$ are diagram-distinct
the sense that $w_{1}$ and $w_{2}$ have different diagrams on the graph $G.$
Also, we could find the necessary condition for the $D_{G}$-freeness of two
arbitrary random variables $a$ and $b.$ i.e., if the supports $\mathbb{F}%
^{+}(G:a)$ and $\mathbb{F}^{+}(G:b)$ are diagram-distinct, in the sense that 
$w_{1}$ and $w_{2}$ are diagram distinct for all pairs $(w_{1},w_{2})$ $\in $
$\mathbb{F}^{+}(G:a)$ $\times $ $\mathbb{F}^{+}(G:b),$ then the random
variables $a$ and $b$ are free over $D_{G}.$
\end{center}

\strut \strut

From Chapter 4 to Chapter 6, we will consider some special $D_{G}$-valued
random variables in a graph $W^{*}$-probability space $\left(
W^{*}(G),E\right) .$ The those random variables are the basic objects to
study Free Probability Theory. We can conclude that

\strut

(i) \ \ if $l$ is a loop, then $L_{l}+L_{l}^{*}$ is $D_{G}$-semicircular.

\strut

(ii) \ if $w$ is a finite path, then $L_{w}+L_{w}^{*}$ is $D_{G}$-even.

\strut

(iii) if $w$ is a finite path, then $L_{w}$ and $L_{w}^{*}$ are $D_{G}$%
-valued R-diagonal.

\strut

In Chapter 5, we consider the generating operator of the graph $W^{*}$%
-algebra $W^{*}(G).$ We compute the moments and cumulants of the generating
operators of the one-vertex graph with $N$-edges and the circulant graph.

\strut \strut

\strut \strut

\strut \strut

\section{Graph $W^{*}$-Probability Spaces}

\strut

\strut

Let $G$ be a countable directed graph and let $\Bbb{F}^{+}(G)$ be the free
semigroupoid of $G.$ i.e., the set $\mathbb{F}^{+}(G)$ is the collection of
all vertices as units and all admissible finite paths of $G.$ Let $w$ be a
finite path with its source $s(w)=x$ and its range $r(w)=y,$ where $x,y\in
V(G).$ Then sometimes we will denote $w$ by $w=xwy$ to express the source
and the range of $w.$ We can define the graph Hilbert space $H_{G}$ by the
Hilbert space $l^{2}\left( \mathbb{F}^{+}(G)\right) $ generated by the
elements in the free semigroupoid $\mathbb{F}^{+}(G).$ i.e., this Hilbert
space has its Hilbert basis $\mathcal{B}=\{\xi _{w}:w\in \mathbb{F}%
^{+}(G)\}. $ Suppose that $w=e_{1}...e_{k}\in FP(G)$ is a finite path with $%
e_{1},...,e_{k}\in E(G).$ Then we can regard $\xi _{w}$ as $\xi
_{e_{1}}\otimes ...\otimes \xi _{e_{k}}.$ So, in [10], Kribs and Power
called this graph Hilbert space the generalized Fock space. Throughout this
paper, we will call $H_{G}$ the graph Hilbert space to emphasize that this
Hilbert space is induced by the graph.

\strut

Define the creation operator $L_{w},$ for $w\in \mathbb{F}^{+}(G),$ by the
multiplication operator by $\xi _{w}$ on $H_{G}.$ Then the creation operator 
$L$ on $H_{G}$ satisfies that

\strut

(i) \ $L_{w}=L_{xwy}=L_{x}L_{w}L_{y},$ for $w=xwy$ with $x,y\in V(G).$

\strut

(ii) $L_{w_{1}}L_{w_{2}}=\left\{ 
\begin{array}{lll}
L_{w_{1}w_{2}} &  & \text{if }w_{1}w_{2}\in \mathbb{F}^{+}(G) \\ 
&  &  \\ 
0 &  & \text{if }w_{1}w_{2}\notin \mathbb{F}^{+}(G),
\end{array}
\right. $

\strut

\ \ \ \ for all $w_{1},w_{2}\in \mathbb{F}^{+}(G).$

\strut

Now, define the annihilation operator $L_{w}^{*},$ for $w\in \mathbb{F}%
^{+}(G)$ by

\strut

\begin{center}
$L_{w}^{\ast }\xi _{w^{\prime }}\overset{def}{=}\left\{ 
\begin{array}{lll}
\xi _{h} &  & \text{if }w^{\prime }=wh\in \mathbb{F}^{+}(G)\xi \\ 
&  &  \\ 
0 &  & \text{otherwise.}
\end{array}
\right. $
\end{center}

\strut

The above definition is gotten by the following observation ;

\strut

\begin{center}
$
\begin{array}{ll}
<L_{w}\xi _{h},\xi _{wh}>\, & =\,<\xi _{wh},\xi _{wh}>\, \\ 
& =\,1=\,<\xi _{h},\xi _{h}> \\ 
& =\,<\xi _{h},L_{w}^{*}\xi _{wh}>,
\end{array}
\,$
\end{center}

\strut

where $<,>$ is the inner product on the graph Hilbert space $H_{G}.$ Of
course, in the above formula we need the admissibility of $w$ and $h$ in $%
\mathbb{F}^{+}(G).$ However, even though $w$ and $h$ are not admissible
(i.e., $wh\notin \mathbb{F}^{+}(G)$), by the definition of $L_{w}^{\ast },$
we have that

\strut

\begin{center}
$
\begin{array}{ll}
<L_{w}\xi _{h},\xi _{h}> & =\,<0,\xi _{h}> \\ 
& =0=\,<\xi _{h},0> \\ 
& =\,<\xi _{h},L_{w}^{*}\xi _{h}>.
\end{array}
\,\,$
\end{center}

\strut

Notice that the creation operator $L$ and the annihilation operator $L^{*}$
satisfy that

\strut

(1.1) \ \ \ $L_{w}^{*}L_{w}=L_{y}$ \ \ and \ \ $L_{w}L_{w}^{*}=L_{x},$ \ for
all \ $w=xwy\in \mathbb{F}^{+}(G),$

\strut

\textbf{under the weak topology}, where $x,y\in V(G).$ Remark that if we
consider the von Neumann algebra $W^{*}(\{L_{w}\})$ generated by $L_{w}$ and 
$L_{w}^{*}$ in $B(H_{G}),$ then the projections $L_{y}$ and $L_{x}$ are
Murray-von Neumann equivalent, because there exists a partial isometry $%
L_{w} $ satisfying the relation (1.1). Indeed, if $w=xwy$ in $\mathbb{F}%
^{+}(G), $ with $x,y\in V(G),$ then under the weak topology we have that

\strut

(1,2) \ \ \ $L_{w}L_{w}^{*}L_{w}=L_{w}$ \ \ and \ \ $%
L_{w}^{*}L_{w}L_{w}^{*}=L_{w}^{*}.$

\strut

So, the creation operator $L_{w}$ is a partial isometry in $W^{*}(\{L_{w}\})$
in $B(H_{G}).$ Assume now that $v\in V(G).$ Then we can regard $v$ as $%
v=vvv. $ So,

\strut

(1.3) $\ \ \ \ \ \ \ \ \ L_{v}^{*}L_{v}=L_{v}=L_{v}L_{v}^{*}=L_{v}^{*}.$

\strut

This relation shows that $L_{v}$ is a projection in $B(H_{G})$ for all $v\in
V(G).$

\strut

Define the \textbf{graph }$W^{*}$\textbf{-algebra} $W^{*}(G)$ by

\strut

\begin{center}
$W^{*}(G)\overset{def}{=}\overline{%
\mathbb{C}[\{L_{w},L_{w}^{*}:w\in
\mathbb{F}^{+}(G)\}]}^{w}.$
\end{center}

\strut

Then all generators are either partial isometries or projections, by (1.2)
and (1.3). So, this graph $W^{\ast }$-algebra contains a rich structure, as
a von Neumann algebra. (This construction can be the generalization of that
of group von Neumann algebra.) Naturally, we can define a von Neumann
subalgebra $D_{G}\subset W^{\ast }(G)$ generated by all projections $L_{v},$ 
$v\in V(G).$ i.e.

\strut

\begin{center}
$D_{G}\overset{def}{=}W^{*}\left( \{L_{v}:v\in V(G)\}\right) .$
\end{center}

\strut

We call this subalgebra the \textbf{diagonal subalgebra} of $W^{*}(G).$
Notice that $D_{G}=\Delta _{\left| G\right| }\subset M_{\left| G\right| }(%
\mathbb{C}),$ where $\Delta _{\left| G\right| }$ is the subalgebra of $%
M_{\left| G\right| }(\mathbb{C})$ generated by all diagonal matrices. Also,
notice that $1_{D_{G}}=\underset{v\in V(G)}{\sum }L_{v}=1_{W^{*}(G)}.$

\strut

If $a\in W^{*}(G)$ is an operator, then it has the following decomposition
which is called the Fourier expansion of $a$ ;

\strut

(1.4) $\ \ \ \ \ \ \ \ \ \ \ a=\underset{w\in \mathbb{F}^{+}(G:a),\,u_{w}\in
\{1,*\}}{\sum }p_{w}L_{w}^{u_{w}},$

\strut

where $p_{w}\in C$ and $\mathbb{F}^{+}(G:a)$ is the support of $a$ defined by

\strut

\begin{center}
$\mathbb{F}^{+}(G:a)=\{w\in \mathbb{F}^{+}(G):p_{w}\neq 0\}.$
\end{center}

\strut

Remark that the free semigroupoid $\mathbb{F}^{+}(G)$ has its partition $%
\{V(G),FP(G)\},$ as a set. i.e.,

\strut

\begin{center}
$\mathbb{F}^{+}(G)=V(G)\cup FP(G)$ \ \ and \ \ $V(G)\cap FP(G)=\emptyset .$
\end{center}

\strut

So, the support of $a$ is also partitioned by

\strut

\begin{center}
$\mathbb{F}^{+}(G:a)=V(G:a)\cup FP(G:a),$
\end{center}

\strut where

\begin{center}
$V(G:a)\overset{def}{=}V(G)\cap \mathbb{F}^{+}(G:a)$
\end{center}

and

\begin{center}
$FP(G:a)\overset{def}{=}FP(G)\cap \mathbb{F}^{+}(G:a).$
\end{center}

\strut

So, the above Fourier expansion (1.4) of the random variable $a$ can be
re-expressed by

\strut

(1.5) $\ \ \ \ \ \ a=\underset{v\in V(G:a)}{\sum }p_{v}L_{v}+\underset{w\in
FP(G:a),\,u_{w}\in \{1,*\}}{\sum }p_{w}L_{w}^{u_{w}}.$

\strut

We can easily see that if $V(G:a)\neq \emptyset ,$ then $\underset{v\in
V(G:a)}{\sum }p_{v}L_{v}$ is contained in the diagonal subalgebra $D_{G}.$
Also, if $V(G:a)=\emptyset ,$ then $\underset{v\in V(G:a)}{\sum }%
p_{v}L_{v}=0_{D_{G}}.$ So, we can define the following canonical conditional
expectation $E:W^{*}(G)\rightarrow D_{G}$ by

\strut

(1.6) \ \ \ $E(a)=E\left( \underset{w\in \mathbb{F}^{+}(G:a),\,u_{w}\in
\{1,*\}}{\sum }p_{w}L_{w}^{u_{w}}\right) \overset{def}{=}\underset{v\in
V(G:a)}{\sum }p_{v}L_{v},$

\strut

for all $a\in W^{*}(G).$ Indeed, $E$ is a well-determined conditional
expectation ; it is a bimodule map satisfying that

\strut

\begin{center}
$E(d)=d,$ for all $d\in D_{G}.$
\end{center}

\strut

And

\strut

\begin{center}
$
\begin{array}{ll}
E\left( dad^{\prime }\right) & =E\left( d(a_{d}+a_{0})d^{\prime }\right)
=E\left( da_{d}d^{\prime }+da_{0}d^{\prime }\right) \\ 
& =E\left( da_{d}d^{\prime }\right) =da_{d}d^{\prime }=d\left( E(a)\right)
d^{\prime },
\end{array}
$
\end{center}

\strut

for all $d,d^{\prime }\in D_{G}$ and $a=a_{d}+a_{0}\in W^{*}(G),$ where

\strut

\begin{center}
$a_{d}=\underset{v\in V(G:a)}{\sum }p_{v}L_{v}$ \ \ and \ \ $a_{0}=\underset{%
w\in FP(G:a),\,u_{w}\in \{1,*\}}{\sum }p_{w}L_{w}^{u_{w}}.$
\end{center}

\strut

Also,

\strut

\begin{center}
$E\left( a^{\ast }\right) =E\left( (a_{d}+a_{0})^{\ast }\right) =E\left(
a_{d}^{\ast }+a_{0}^{\ast }\right) =a_{d}^{\ast }=E(a)^{\ast },$
\end{center}

\strut

for all $a\in W^{*}(G).$ Here, $a_{d}^{*}=\left( \underset{v\in V(G:a)}{\sum 
}p_{v}L_{v}\right) ^{*}=\underset{v\in V(G:a)}{\sum }\overline{p_{v}}\,L_{v}$
in $D_{G}.$

\strut \strut \strut \strut

\begin{definition}
Let $G$ be a countable directed graph and let $W^{*}(G)$ be the graph $W^{*}$%
-algebra induced by $G.$ Let $E:W^{*}(G)\rightarrow D_{G}$ be the
conditional expectation defined above. Then we say that the algebraic pair $%
\left( W^{*}(G),E\right) $ is the graph $W^{*}$-probability space over the
diagonal subalgebra $D_{G}$. By the very definition, it is one of the $W^{*}$%
-probability space with amalgamation over $D_{G}.$ All elements in $\left(
W^{*}(G),E\right) $ are called $D_{G}$-valued random variables.
\end{definition}

\strut

We have a graph $W^{*}$-probability space $\left( W^{*}(G),E\right) $ over
its diagonal subalgebra $D_{G}.$ We will define the following free
probability data of $D_{G}$-valued random variables.

\strut

\begin{definition}
Let $W^{*}(G)$ be the graph $W^{*}$-algebra induced by $G$ and let $a\in
W^{*}(G).$ Define the $n$-th ($D_{G}$-valued) moment of $a$ by

\strut 

$\ \ \ \ \ E\left( d_{1}ad_{2}a...d_{n}a\right) ,$ for all $n\in \mathbb{N}$,

\strut 

where $d_{1},...,d_{n}\in D_{G}$. Also, define the $n$-th ($D_{G}$-valued)
cumulant of $a$ by

\strut 

$\ \ \ \ \ k_{n}(d_{1}a,d_{2}a,...,d_{n}a)=C^{(n)}\left( d_{1}a\otimes
d_{2}a\otimes ...\otimes d_{n}a\right) ,$

\strut 

for all $n\in \mathbb{N},$ and for $d_{1},...,d_{n}\in D_{G},$ where $%
\widehat{C}=(C^{(n)})_{n=1}^{\infty }\in I^{c}\left( W^{*}(G),D_{G}\right) $
is the cumulant multiplicative bimodule map induced by the conditional
expectation $E,$ in the sense of Speicher. We define the $n$-th trivial
moment of $a$ and the $n$-th trivial cumulant of $a$ by

\strut 

$\ \ \ \ \ E(a^{n})$ $\ \ $and $\ \ k_{n}\left( \underset{n-times}{%
\underbrace{a,a,...,a}}\right) =C^{(n)}\left( a\otimes a\otimes ...\otimes
a\right) ,$

\strut 

respectively, for all $n\in \mathbb{N}.$
\end{definition}

\strut

To compute the $D_{G}$-valued moments and cumulants of the $D_{G}$-valued
random variable $a,$ we need to introduce the following new definition ;

\strut

\begin{definition}
Let $\left( W^{*}(G),E\right) $ be a graph $W^{*}$-probability space over $%
D_{G}$ and let $a\in \left( W^{*}(G),E\right) $ be a random variable. Define
the subset $FP_{*}(G:a)$ in $FP(G:a)$ \ by

\strut 

$\ \ \ FP_{*}\left( G:a\right) \overset{def}{=}\{w\in \mathbb{F}^{+}(G:a):$%
both $L_{w}$ and $L_{w}^{*}$ are summands of $a\}.$

\strut 

And let $FP_{*}^{c}(G:a)\overset{def}{=}FP(G:a)\,\setminus \,FP_{*}(G:a).$
\end{definition}

\strut \strut \strut

We already observed that if $a\in \left( W^{*}(G),E\right) $ is a $D_{G}$%
-valued random variable, then $a$ has its Fourier expansion $a_{d}+a_{0},$
where

\strut

\begin{center}
$a_{d}=\underset{v\in V(G:a)}{\sum }p_{v}L_{v}$
\end{center}

and

\begin{center}
$a_{0}=\underset{w\in FP(G:a),\,u_{w}\in \{1,*\}}{\sum }p_{w}L_{w}^{u_{w}}.$
\end{center}

\strut

By the previous definition, the set $FP(G:a)$ is partitioned by

\strut

\begin{center}
$FP(G:a)=FP_{*}(G:a)\cup FP_{*}^{c}(G:a),$
\end{center}

\strut

for the fixed random variable $a$ in $\left( W^{*}(G),E\right) .$ So, the
summand $a_{0},$ in the Fourier expansion of $a=a_{d}+a_{0},$ has the
following decomposition ;

\strut

\begin{center}
$a_{0}=a_{(*)}+a_{(non-*)},$
\end{center}

\strut where\strut

\begin{center}
$a_{(*)}=\underset{l\in FP_{*}(G:a)}{\sum }\left(
p_{l}L_{l}+p_{l^{t}}L_{l}^{*}\right) $
\end{center}

and

\begin{center}
$a_{(non-*)}=\underset{w\in FP_{*}^{c}(G:a),\,u_{w}\in \{1,*\}}{\sum }%
p_{w}L_{w}^{u_{w}},$
\end{center}

\strut

where $p_{l^{t}}$ is the coefficient of $L_{l}^{\ast }$ depending on $l\in
FP_{\ast }(G:a).$ (There is no special meaning for the complex number $%
p_{l^{t}}.$ But we have to keep in mind that $p_{l}\neq p_{l^{t}},$ in
general. i.e. $a_{(\ast )}=\underset{l_{1}\in FP_{\ast }(G:a)}{\sum }%
p_{l_{1}}L_{l_{1}}+\underset{l_{2}\in FP_{\ast }(G:a)}{\sum }%
p_{l_{2}}L_{l_{2}}^{\ast }$ ! But for the convenience of using notation, we
will use the notation $p_{l^{t}},$ for the coefficient of $L_{l}^{\ast }.$)
For instance, let $V(G:a)=\{v_{1},v_{2}\}$ and $FP(G:a)=\{w_{1},w_{2}\}$ and
let the random variable $a$ in $\left( W^{\ast }(G),E\right) $ be

\strut

\begin{center}
$a=L_{v_{1}}+L_{v_{2}}+L_{w_{1}}^{*}+L_{w_{1}}+L_{w_{2}}^{*}.$
\end{center}

\strut

Then we have that $a_{d}=L_{v_{1}}+L_{v_{2}}$, $%
a_{(*)}=L_{w_{1}}^{*}+L_{w_{1}}$ and $a_{(non-*)}=L_{w_{2}}^{*}.$ By
definition, $a_{0}=a_{(*)}+a_{(non-*)}.$

\strut \strut \strut

\strut \strut

\strut \strut

\section{$D_{G}$-Moments and $D_{G}$-Cumulants of Random Variables}

\strut

\strut

\strut

Throughout this chapter, let $G$ be a countable directed graph and let $%
\left( W^{*}(G),E\right) $ be the graph $W^{*}$-probability space over its
diagonal subalgebra $D_{G}.$ In this chapter, we will compute the $D_{G}$%
-valued moments and the $D_{G}$-valued cumulants of arbitrary random variable

$\strut $

\begin{center}
$a=\underset{w\in \mathbb{F}^{+}(G:a),\,u_{w}\in \{1,*\}}{\sum }%
p_{w}L_{w}^{u_{w}}$
\end{center}

\strut

in the graph $W^{*}$-probability space $\left( W^{*}(G),E\right) $.

\strut

\strut

\strut

\subsection{Lattice Path Model}

\strut

\strut

\strut

Throughout this section, let $G$ be a countable directed graph and let $%
\left( W^{*}(G),E\right) $ be the graph $W^{*}$-probability space over its
diagonal subalgebra $D_{G}.$ Let $w_{1},...,w_{n}\in \Bbb{F}^{+}(G)$ and let 
$L_{w_{1}}^{u_{w_{1}}}...L_{w_{n}}^{u_{w_{n}}}\in \left( W^{*}(G),E\right) $
be a $D_{G}$-valued random variable. In this section, we will define a
lattice path model for the random variable $%
L_{w_{1}}^{u_{w_{1}}}...L_{w_{n}}^{u_{w_{n}}}.$ Recall that if $%
w=e_{1}....e_{k}\in FP(G)$ with $e_{1},...,e_{k}\in E(G),$ then we can
define the length $\left| w\right| $ of $w$ by $k.$ i.e.e, the length $%
\left| w\right| $ of $w$ is the cardinality $k$ of the admissible edges $%
e_{1},...,e_{k}.$

\strut

\begin{definition}
Let $G$ be a countable directed graph and $\Bbb{F}^{+}(G),$ the free
semigroupoid. If $w\in \Bbb{F}^{+}(G),$ then $L_{w}$ is the corresponding $%
D_{G}$-valued random variable in $\left( W^{*}(G),E\right) .$ We define the
lattice path $l_{w}$ of $L_{w}$ and the lattice path $l_{w}^{-1}$ of $%
L_{w}^{*}$ by the lattice paths satisfying that ;

\strut 

(i) \ \ the lattice path $l_{w}$ starts from $*=(0,0)$ on the $\Bbb{R}^{2}$%
-plane.

\strut 

(ii) \ if $w\in V(G),$ then $l_{w}$ has its end point $(0,1).$

\strut 

(iii) if $w\in E(G),$ then $l_{w}$ has its end point $(1,1).$

\strut 

(iv) if $w\in E(G),$ then $l_{w}^{-1}$ has its end point $(-1,-1).$

\strut 

(v) \ if $w\in FP(G)$ with $\left| w\right| =k,$ then $l_{w}$ has its end
point $(k,k).$

\strut 

(vi) if $w\in FP(G)$ with $\left| w\right| =k,$ then $l_{w}^{-1}$ has its
end point $(-k,-k).$

\strut 

Assume that finite paths $w_{1},...,w_{s}$ in $FP(G)$ satisfy that $%
w_{1}...w_{s}\in FP(G).$ Define the lattice path $l_{w_{1}...w_{s}}$ by the
connected lattice path of the lattice paths $l_{w_{1}},$ ..., $l_{w_{s}}.$
i.e.e, $l_{w_{2}}$ starts from $(k_{w_{1}},k_{w_{1}})\in \Bbb{R}^{+}$ and
ends at $(k_{w_{1}}+k_{w_{2}},k_{w_{1}}+k_{w_{2}}),$ where $\left|
w_{1}\right| =k_{w_{1}}$ and $\left| w_{2}\right| =k_{w_{2}}.$ Similarly, we
can define the lattice path $l_{w_{1}...w_{s}}^{-1}$ as the connected path
of $l_{w_{s}}^{-1},$ $l_{w_{s-1}}^{-1},$ ..., $l_{w_{1}}^{-1}.$
\end{definition}

\strut

\begin{definition}
Let $G$ be a countable directed graph and assume that $%
L_{w_{1}},...,L_{w_{n}}$ are generators of $\left( W^{*}(G),E\right) .$ Then
we have the lattice paths $l_{w_{1}},$ ..., $l_{w_{n}}$ of $L_{w_{1}},$ ..., 
$L_{w_{n}},$ respectively in $\Bbb{R}^{2}.$ Suppose that $%
L_{w_{1}}^{u_{w_{1}}}...L_{w_{n}}^{u_{w_{n}}}\neq 0_{D_{G}}$ in $\left(
W^{*}(G),E\right) ,$ where $u_{w_{1}},...,u_{w_{n}}\in \{1,*\}.$ Define the
lattice path $l_{w_{1},...,w_{n}}^{u_{w_{1}},...,u_{w_{n}}}$ of nonzero $%
L_{w_{1}}^{u_{w_{1}}}...L_{w_{n}}^{u_{w_{n}}}$ by the connected lattice path
of $l_{w_{1}}^{t_{w_{1}}},$ ..., $l_{w_{n}}^{t_{w_{n}}},$ where $t_{w_{j}}=1$
if $u_{w_{j}}=1$ and $t_{w_{j}}=-1$ if $u_{w_{j}}=*.$ Assume that $%
L_{w_{1}}^{u_{w_{1}}}...L_{w_{n}}^{u_{w_{n}}}$ $=$ $0_{D_{G}}.$ Then the
empty set $\emptyset $ in $\Bbb{R}^{2}$ is the lattice path of it. We call
it the empty lattice path. By $LP_{n},$ we will denote the set of all
lattice paths of the $D_{G}$-valued random variables having their forms of $%
L_{w_{1}}^{u_{w_{1}}}...L_{w_{n}}^{u_{w_{n}}},$ including empty lattice path.
\end{definition}

\strut

Also, we will define the following important property on the set of all
lattice paths ;

\strut

\begin{definition}
Let $l_{w_{1},...,w_{n}}^{u_{w_{1}},...,u_{w_{n}}}\neq \emptyset $ be a
lattice path of $L_{w_{1}}^{u_{w_{1}}}...L_{w_{n}}^{u_{w_{n}}}\neq 0_{D_{G}}$
in $LP_{n}.$ If the lattice path $%
l_{w_{1},...,w_{n}}^{u_{w_{1}},...,u_{w_{n}}}$ starts from $*$ and ends on
the $*$-axis in $\Bbb{R}^{+},$ then we say that the lattice path $%
l_{w_{1},...,w_{n}}^{u_{w_{1}},...,u_{w_{n}}}$ has the $*$-axis-property. By 
$LP_{n}^{*},$ we will denote the set of all lattice paths having their forms
of $l_{w_{1},...,w_{n}}^{u_{w_{1}},...,u_{w_{n}}}$ which have the $*$%
-axis-property. By little abuse of notation, sometimes, we will say that the 
$D_{G}$-valued random variable $L_{w_{1}}^{u_{w_{1}}}...L_{w_{n}}^{u_{w_{n}}}
$satisfies the $*$-axis-property if the lattice path $%
l_{w_{1},...,w_{n}}^{u_{w_{1}},...,u_{w_{n}}}$ of it has the $*$%
-axis-property.
\end{definition}

\strut

The following theorem shows that finding $E\left(
L_{w_{1}}^{u_{w_{1}}}...L_{w_{n}}^{u_{w_{n}}}\right) $ is checking the $*$%
-axis-property of $L_{w_{1}}^{u_{w_{1}}}...L_{w_{n}}^{u_{w_{n}}}.$

\strut

\begin{theorem}
Let $L_{w_{1}}^{u_{w_{1}}}...L_{w_{n}}^{u_{w_{n}}}\in \left(
W^{*}(G),E\right) $ be a $D_{G}$-valued random variable, where $%
u_{w_{1}},...,u_{w_{n}}\in \{1,*\}.$ Then $E\left(
L_{w_{1}}^{u_{w_{1}}}...L_{w_{n}}^{u_{w_{n}}}\right) $ $\neq $ $0_{D_{G}}$
if and only if $L_{w_{1}}^{u_{w_{1}}}...L_{w_{n}}^{u_{w_{n}}}$ has the $*$%
-axis-property (i.e., the corresponding lattice path $%
l_{w_{1},...,w_{n}}^{u_{w_{1}},...,u_{w_{n}}}$ of $%
L_{w_{1}}^{u_{w_{1}}}...L_{w_{n}}^{u_{w_{n}}}$ is contained in $LP_{n}^{*}.$
Notice that $\emptyset \notin LP_{n}^{*}.$)
\end{theorem}

\strut

\begin{proof}
($\Leftarrow $) Let $l=l_{w_{1},...,w_{n}}^{u_{w_{1}},...,u_{w_{n}}}\in
LP_{n}^{*}.$ Suppose that $w_{1}=vw_{1}v_{1}^{\prime }$ and $%
w_{n}=v_{n}w_{n}v_{n}^{\prime },$ for $v_{1},$ $v_{1}^{\prime },$ $v_{n},$ $%
v_{n}^{\prime }$ $\in $ $V(G).$ If $l$ is in $LP_{n}^{*},$ then

\strut

(2.1.1)$\ \ \ \ \ \ \ \left\{ 
\begin{array}{lll}
v_{1}=v_{n}^{\prime } &  & \text{if }u_{w_{1}}=1\text{ and }u_{w_{n}}=1 \\ 
&  &  \\ 
v_{1}=v_{n} &  & \text{if }u_{w_{1}}=1\text{ and }u_{w_{n}}=* \\ 
&  &  \\ 
v_{1}^{\prime }=v_{n}^{\prime } &  & \text{if }u_{w_{1}}=*\text{ and }%
u_{w_{n}}=1 \\ 
&  &  \\ 
v_{1}^{\prime }=v_{n} &  & \text{if }u_{w_{1}}=*\text{ and }u_{w_{n}}=*.
\end{array}
\right. $

\strut

By the definition of lattice paths having the $*$-axis-property and by
(2.1.1), if $l_{w_{1},...,w_{n}}^{u_{w_{1}},...,u_{w_{n}}}\in LP_{n}^{*},$
then there exists $v\in V(G)$ such that

\strut

$\ \ \ \ \ \ \ \ \ \ \ \ \ \ \ \ \ \ \
L_{w_{1}}^{u_{w_{1}}}...L_{w_{n}}^{u_{w_{n}}}=L_{v},$

where

(2.1.2) \ \ \ \ $\ \ \left\{ 
\begin{array}{ll}
v=v_{1}=v_{n}^{\prime } & \text{if }u_{w_{1}}=1\text{ and }u_{w_{n}}=1 \\ 
&  \\ 
v=v_{1}=v_{n} & \text{if }u_{w_{1}}=1\text{ and }u_{w_{n}}=* \\ 
&  \\ 
v=v_{1}^{\prime }=v_{n}^{\prime } & \text{if }u_{w_{1}}=*\text{ and }%
u_{w_{n}}=1 \\ 
&  \\ 
v=v_{1}^{\prime }=v_{n} & \text{if }u_{w_{1}}=*\text{ and }u_{w_{n}}=*.
\end{array}
\right. $

\strut

This shows that $E\left(
L_{w_{1}}^{u_{w_{1}}}...L_{w_{n}}^{u_{w_{n}}}\right) =L_{v}\neq 0_{D_{G}}.$

\strut

($\Rightarrow $) Assume that $E\left(
L_{w_{1}}^{u_{w_{1}}}...L_{w_{n}}^{u_{w_{n}}}\right) \neq 0_{D_{G}}.$ This
means that there exists $L_{v},$ with $v\in V(G),$ such that

\strut

(2.1.3) \ \ \ \ \ \ \ \ $\ \ \ E\left(
L_{w_{1}}^{u_{w_{1}}}...L_{w_{n}}^{u_{w_{n}}}\right) =L_{v}.$

\strut

Equivalently, we have that $%
L_{w_{1}}^{u_{w_{1}}}...L_{w_{n}}^{u_{w_{n}}}=L_{v}$ in $W^{*}(G).$ Let $%
l=l_{w_{1},...,w_{n}}^{u_{w_{1}},...,u_{w_{n}}}\in LP_{n}$ be the lattice
path of the $D_{G}$-valued random variable $%
L_{w_{1}}^{u_{w_{1}}}...L_{w_{n}}^{u_{w_{n}}}.$ By (2.1.3), trivially, $%
l\neq \emptyset ,$ since $l$ should be the connected lattice path. Assume
that the nonempty lattice path $l$ is contained in $LP_{n}\,\,\setminus
\,LP_{n}^{*}.$ Then, under the same conditions of (2.1.1), we have that

\strut

(2.1.4) \ \ \ \ \ \ $\left\{ 
\begin{array}{lll}
v_{1}\neq v_{n}^{\prime } &  & \text{if }u_{w_{1}}=1\text{ and }u_{w_{n}}=1
\\ 
&  &  \\ 
v_{1}\neq v_{n} &  & \text{if }u_{w_{1}}=1\text{ and }u_{w_{n}}=* \\ 
&  &  \\ 
v_{1}^{\prime }\neq v_{n}^{\prime } &  & \text{if }u_{w_{1}}=*\text{ and }%
u_{w_{n}}=1 \\ 
&  &  \\ 
v_{1}^{\prime }\neq v_{n} &  & \text{if }u_{w_{1}}=*\text{ and }u_{w_{n}}=*.
\end{array}
\right. $

\strut

Therefore, by (2.1.2), there is no vertex $v$ satisfying $%
L_{w_{1}}^{u_{w_{1}}}...L_{w_{n}}^{u_{w_{n}}}=L_{v}.$ This contradict our
assumption.
\end{proof}

\strut

By the previous theorem, we can conclude that $E\left(
L_{w_{1}}^{u_{w_{1}}}...L_{w_{n}}^{u_{w_{n}}}\right) =L_{v},$ for some $v\in
V(G)$ if and only if the lattice path $%
l_{w_{1},...,w_{n}}^{u_{w_{1}},...,u_{w_{n}}}$ has the $*$-axis-property
(i.e., $l_{w_{1},...,w_{n}}^{u_{w_{1}},...,u_{w_{n}}}\in LP_{n}^{*}$).\strut
\strut \strut

\strut \strut

\subsection{$D_{G}$-Valued Moments and Cumulants of Random Variables\strut}

\bigskip

\bigskip

Let $w_{1},...,w_{n}\in \mathbb{F}^{+}(G)$, $u_{1},...,u_{n}\in \{1,*\}$ and
let $L_{w_{1}}^{u_{1}}...L_{w_{n}}^{u_{n}}\in \left( W^{*}(G),E\right) $ be
a $D_{G}$-valued random variable. Recall that, in the previous section, we
observed that the $D_{G}$-valued random variable $%
L_{w_{1}}^{u_{1}}...L_{w_{n}}^{u_{n}}=L_{v}\in \left( W^{*}(G),E\right) $
with $v\in V(G)$ if and only if the lattice path $%
l_{w_{1},...,w_{n}}^{u_{1},...,u_{n}}$ of \ $%
L_{w_{1}}^{u_{1}}...L_{w_{n}}^{u_{n}}$ has the $*$-axis-property
(equivalently, $l_{w_{1},...,w_{n}}^{u_{1},...,u_{n}}\in LP_{n}^{*}$).
Throughout this section, fix a $D_{G}$-valued random variable $a\in \left(
W^{*}(G),E\right) .$ Then the $D_{G}$-valued random variable $a$ has the
following Fourier expansion,

\bigskip

\begin{center}
$a=\underset{v\in V(G:a)}{\sum }p_{v}L_{v}+\underset{l\in FP_{*}(G:a)}{\sum }%
\left( p_{l}L_{l}+p_{l^{t}}L_{l}\right) +\underset{w\in
FP_{*}^{c}(G:a),~u_{w}\in \{1,*\}}{\sum }p_{w}L_{w}^{u_{w}}.$
\end{center}

\bigskip

Let's observe the new $D_{G}$-valued random variable $d_{1}ad_{2}a...d_{n}a%
\in \left( W^{*}(G),E\right) ,$ where $d_{1},...,d_{n}\in D_{G}$ and $a\in
W^{*}(G)$ is given. Put

\strut

\begin{center}
$d_{j}=\underset{v_{j}\in V(G:d_{j})}{\sum }q_{v_{j}}L_{v_{j}}\in D_{G},$
for \ $j=1,...,n.$
\end{center}

\strut

Notice that $V(G:d_{j})=\mathbb{F}^{+}(G:d_{j}),$ since $d_{j}\in
D_{G}\hookrightarrow W^{\ast }(G).$ Then

\strut

$\ d_{1}ad_{2}a...d_{n}a$

\strut

$\ \ \ =\left( \underset{v_{1}\in V(G:d_{1})}{\sum }q_{v_{1}}L_{v_{1}}%
\right) \left( \underset{w_{1}\in \mathbb{F}^{+}(G:a),\,u_{w_{1}}\in \{1,*\}%
}{\sum }p_{w_{1}}L_{w_{1}}^{u_{w_{1}}}\right) $

$\ \ \ \ \ \ \ \ \ \ \ \ \ \ \ \ \ \cdot \cdot \cdot \left( \underset{%
v_{1}\in V(G:d_{n})}{\sum }q_{v_{n}}L_{v_{n}}\right) \left( \underset{%
w_{n}\in \mathbb{F}^{+}(G:a),\,u_{w_{n}}\in \{1,*\}}{\sum }%
p_{w_{n}}L_{w_{n}}^{u_{w_{n}}}\right) $

\strut

$\ \ \ =\underset{(v_{1},...,v_{n})\in \Pi _{j=1}^{n}V(G:d_{j})}{\sum }%
\left( q_{v_{1}}...q_{v_{n}}\right) $

$\ \ \ \ \ \ \ \ \ \ \ \ \ \ \ \ \ \ \ \ \ \ \ \ \ \ \ \ \ \ \ \
(L_{v_{1}}\left( \underset{w_{1}\in \mathbb{F}^{+}(G:a),\,u_{w_{1}}\in
\{1,*\}}{\sum }p_{w_{1}}L_{w_{1}}^{u_{w_{1}}}\right) $

$\ \ \ \ \ \ \ \ \ \ \ \ \ \ \ \ \ \ \ \ \ \ \ \ \ \ \ \ \ \ \ \ \ \ \ \ \ \
\ \cdot \cdot \cdot L_{v_{n}}\left( \underset{w_{n}\in \mathbb{F}%
^{+}(G:a),\,u_{w_{n}}\in \{1,*\}}{\sum }p_{w_{n}}L_{w_{n}}^{u_{w_{n}}}%
\right) )$

\strut

(1.2.1)

\strut

$\ \ \ =\underset{(v_{1},...,v_{n})\in \Pi _{j=1}^{n}V(G:d_{j})}{\sum }%
\left( q_{v_{1}}...q_{v_{n}}\right) $

\strut

$\ \ \ \ \ \underset{(w_{1},...,w_{n})\in \mathbb{F}^{+}(G:a)^{n},\,%
\,u_{w_{j}}\in \{1,*\}}{\sum }\left( p_{w_{1}}...p_{w_{n}}\right)
L_{v_{1}}L_{w_{1}}^{u_{w_{1}}}...L_{v_{n}}L_{w_{n}}^{u_{w_{n}}}.$

\strut

Now, consider the random variable $%
L_{v_{1}}L_{w_{1}}^{u_{w_{1}}}...L_{v_{n}}L_{w_{n}}^{u_{w_{n}}}$ in the
formula (1.2.1). Suppose that $w_{j}=x_{j}w_{j}y_{j},$ with $x_{j},y_{j}\in
V(G),$ for all \ $j=1,...,n.$ Then

\strut

\strut (1.2.2)

\begin{center}
$L_{v_{1}}L_{w_{1}}^{u_{w_{1}}}...L_{v_{n}}L_{w_{n}}^{u_{w_{n}}}=\delta
_{(v_{1},x_{1},y_{1}:u_{w_{1}})}\cdot \cdot \cdot \delta
_{(v_{n},x_{n},y_{n}:u_{w_{n}})}\left(
L_{w_{1}}^{u_{w_{n}}}...L_{w_{n}}^{u_{w_{n}}}\right) ,$
\end{center}

\strut

where

\strut

\begin{center}
$\delta _{(v_{j},x_{j},y_{j}:u_{w_{j}})}=\left\{ 
\begin{array}{lll}
\delta _{v_{j},x_{j}} &  & \text{if }u_{w_{j}}=1 \\ 
&  &  \\ 
\delta _{v_{j},y_{j}} &  & \text{if }u_{w_{j}}=*,
\end{array}
\right. $
\end{center}

\strut

for all \ $j=1,...,n,$ where $\delta $ in the right-hand side is the
Kronecker delta. So, the left-hand side can be understood as a (conditional)
Kronecker delta depending on $\{1,*\}$.

\strut \strut

By (1.2.1) and (1.2.2), the $n$-th moment of $a$ is

\strut

$\ E\left( d_{1}a...d_{n}a\right) $

\strut

\ \ $\ \ \ =E(\underset{(v_{1},...,v_{n})\in \Pi _{j=1}^{n}V(G:d_{j})}{\sum }%
\left( \Pi _{j=1}^{n}q_{v_{j}}\right) $

\strut

$\ \ \ \ \ \ \ \ \ \ \ \ \ \ \ \underset{(w_{1},...,w_{n})\in \mathbb{F}%
^{+}(G:a)^{n},\,w_{j}=x_{j}w_{j}y_{j},\,u_{w_{j}}\in \{1,\ast \}}{\sum }%
\left( \Pi _{j=1}^{n}p_{w_{j}}\right) $

\strut

$\ \ \ \ \ \ \ \ \ \ \ \ \ \ \ \ \ \ \ \ \ \ \ \ \ \ \ \ \ \ \ \ \ \ \left(
\Pi _{j=1}^{n}\delta _{(v_{j},x_{j},y_{j}:u_{w_{j}})}\right) \left(
L_{w_{1}}^{u_{w_{1}}}...L_{w_{n}}^{u_{w_{n}}}\right) )$

\strut

$\ \ \ \ \ \ =\underset{(v_{1},...,v_{n})\in \Pi _{j=1}^{n}V(G:d_{j})}{\sum }%
\left( \Pi _{j=1}^{n}q_{v_{j}}\right) $

\strut

$\ \ \ \ \ \ \ \ \ \ \underset{(w_{1},...,w_{n})\in \mathbb{F}%
^{+}(G:a)^{n},\,w_{j}=x_{j}w_{j}y_{j},\,u_{w_{j}}\in \{1,\ast \}}{\sum }%
\left( \Pi _{j=1}^{n}p_{w_{j}}\right) $

\strut

$\ \ \ \ \ \ \ \ \ \ \ \ \ \ \ \ \ \ \ \ \ \ \ \ \ \ \ \ \ \ \ \ \left( \Pi
_{j=1}^{n}\delta _{(v_{j},x_{j},y_{j}:u_{w_{j}})}\right) \,\,E\left(
L_{w_{1}}^{u_{w_{1}}}...L_{w_{n}}^{u_{w_{n}}}\right) .$

\strut

Thus to compute the $n$-th moment of $a$, we have to observe $E\left(
L_{w_{1}}^{u_{w_{1}}}...L_{w_{n}}^{u_{w_{n}}}\right) .$ In the previous
section, we observed that $E\left(
L_{w_{1}}^{u_{w_{1}}}...L_{w_{n}}^{u_{w_{n}}}\right) $ is nonvanishing if
and only if $L_{w_{1}}^{u_{w_{1}}}...L_{w_{n}}^{u_{w_{n}}}$ has the $*$%
-axis-property.

\strut

\begin{proposition}
Let $a\in \left( W^{*}(G),E\right) $ be given as above. Then the $n$-th
moment of $a$ is

\strut 

$\ \ E\left( d_{1}a...d_{n}a\right) =\underset{(v_{1},...,v_{n})\in \Pi
_{j=1}^{n}V(G:d_{j})}{\sum }\left( \Pi _{j=1}^{n}q_{v_{j}}\right) $

\strut 

$\ \ \ \underset{(w_{1},...,w_{n})\in \mathbb{F}^{+}(G:a)^{n},\,u_{w_{j}}\in
\{1,*\},\,l_{w_{1},...,w_{n}}^{u_{w_{1}},...,u_{w_{n}}}\in LP_{n}^{*}}{\sum }%
\left( \Pi _{j=1}^{n}p_{w_{j}}\right) $

\strut 

$\ \ \ \ \ \ \ \ \ \ \ \ \ \ \ \left( \Pi _{j=1}^{n}\delta
_{(v_{j},x_{j},y_{j}:u_{w_{j}})}\right) \,\,E\left(
L_{w_{1}}^{u_{w_{1}}}...L_{w_{n}}^{u_{w_{n}}}\right) .$

$\square $
\end{proposition}

\strut

From now, rest of this paper, we will compute the $D_{G}$-valued cumulants
of the given $D_{G}$-valued random variable $a$. Let $w_{1},...,w_{n}\in
FP(G)$ be finite paths and \ $u_{1},...,u_{n}\in \{1,*\}$. Then, by the
M\"{o}bius inversion, we have

\strut

(2.2.1)$\ \ $

\begin{center}
$k_{n}\left( L_{w_{1}}^{u_{1}}~,...,~L_{w_{n}}^{u_{n}}\right) =\underset{\pi
\in NC(n)}{\sum }\widehat{E}(\pi )\left( L_{w_{1}}^{u_{1}}~\otimes
...\otimes ~L_{w_{n}}^{u_{n}}\right) \mu (\pi ,1_{n}),$
\end{center}

\strut

where $\widehat{E}=\left( E^{(n)}\right) _{n=1}^{\infty }$ is the moment
multiplicative bimodule map induced by the conditional expectation $E$ (See
[16]) and where $NC(n)$ is the collection of all noncrossing partition over $%
\{1,...,n\}.$ Notice that if $L_{w_{1}}^{u_{1}}...L_{w_{n}}^{u_{n}}$ does
not have the $*$-axis-property, then

\strut

\begin{center}
$E\left( L_{w_{1}}^{u_{1}}...L_{w_{n}}^{u_{n}}\right) =0_{D_{G}},$
\end{center}

\strut

by Section 2.1. Consider the noncrossing partition $\pi \in NC(n)$ with its
blocks $V_{1},...,V_{k}$. Choose one block $V_{j}=(j_{1},...,j_{k})\in \pi .$
Then we have that

\strut

(2.2.2) $\ $

\begin{center}
$\widehat{E}(\pi \mid _{V_{j}})\left( L_{w_{1}}^{u_{1}}~\otimes ...\otimes
~L_{w_{n}}^{u_{n}}\right) =E\left(
L_{w_{j_{1}}}^{u_{j_{1}}}d_{j_{2}}L_{w_{j_{2}}}^{u_{j_{2}}}...d_{j_{k}}L_{w_{j_{k}}}^{u_{j_{k}}}\right) , 
$
\end{center}

\strut where

\begin{center}
$d_{j_{i}}=\left\{ 
\begin{array}{ll}
1_{D_{G}} & 
\begin{array}{l}
\text{if there is no inner blocks} \\ 
\text{between }j_{i-1}\text{ \ and }j_{i}\text{ in }V_{j}
\end{array}
\\ 
&  \\ 
L_{v_{j_{i}}}\neq 1_{D_{G}} & 
\begin{array}{l}
\text{if there are inner blocks} \\ 
\text{between }j_{i-1}\text{ \ and }j_{i}\text{ in }V_{j},
\end{array}
\end{array}
\right. $
\end{center}

\strut

where $v_{j_{2}},...,v_{j_{k}}\in V(G).$ So, again by Section 2.1, $\widehat{%
E}(\pi \mid _{V_{j}})\left( L_{w_{1}}^{u_{1}}~\otimes ...\otimes
~L_{w_{n}}^{u_{n}}\right) $ is nonvanishing if and only if $%
L_{w_{j_{1}}}^{u_{j_{1}}}d_{j_{2}}L_{w_{j_{2}}}^{u_{j_{2}}}...d_{j_{k}}L_{w_{j_{k}}}^{u_{j_{k}}} 
$ has the $*$-axis-property, for all \ $j=1,...,n.$

\strut

Assume that

\begin{center}
$\widehat{E}(\pi \mid _{V_{j}})\left( L_{w_{1}}^{u_{1}}~\otimes ...\otimes
~L_{w_{n}}^{u_{n}}\right) =L_{v_{j}}$
\end{center}

and

\begin{center}
$\widehat{E}(\pi \mid _{V_{i}})\left( L_{w_{1}}^{u_{1}}~\otimes ...\otimes
~L_{w_{n}}^{u_{n}}\right) =L_{v_{i}}.$
\end{center}

\strut

If $v_{j}\neq v_{i},$ then the partition-dependent $D_{G}$-moment satisfies
that

\strut

\begin{center}
$\widehat{E}(\pi )\left( L_{w_{1}}^{u_{1}}~\otimes ...\otimes
~L_{w_{n}}^{u_{n}}\right) =0_{D_{G}}.$
\end{center}

\strut

This says that $\widehat{E}(\pi )\left( L_{w_{1}}^{u_{1}}~\otimes ...\otimes
~L_{w_{n}}^{u_{n}}\right) \neq 0_{D_{G}}$ if and only if there exists $v\in
V(G)$ such that

\strut

\begin{center}
$\widehat{E}(\pi \mid _{V_{j}})\left( L_{w_{1}}^{u_{1}}~\otimes ...\otimes
~L_{w_{n}}^{u_{n}}\right) =L_{v},$
\end{center}

\strut

for all \ $j=1,...,k.$

\strut

\begin{definition}
Let $NC(n)$ be the set of all \ noncrossing partition over $\{1,...,n\}$ and
let $L_{w_{1}}^{u_{1}},$ $...,$ $L_{w_{n}}^{u_{n}}\in \left(
W^{*}(G),E\right) $ be $D_{G}$-valued random variables, where $%
u_{1},...,u_{n}\in \{1,*\}.$ We say that the $D_{G}$-valued random variable $%
L_{w_{1}}^{u_{1}}...L_{w_{n}}^{u_{n}}$ is $\pi $-connected if the $\pi $%
-dependent $D_{G}$-moment of it is nonvanishing, for $\pi \in NC(n).$ In
other words, the random variable $L_{w_{1}}^{u_{1}}...L_{w_{n}}^{u_{n}}$ is $%
\pi $-connected, for $\pi \in NC(n),$ if

\strut 

$\ \ \ \ \ \ \ \ \ \ \ \widehat{E}(\pi )\left( L_{w_{1}}^{u_{1}}~\otimes
...\otimes ~L_{w_{n}}^{u_{n}}\right) \neq 0_{D_{G}}.$

\strut 

i.e., there exists a vertex $v\in V(G)$ such that

\strut 

$\ \ \ \ \ \ \ \ \ \ \ \widehat{E}(\pi )\left( L_{w_{1}}^{u_{1}}~\otimes
...\otimes ~L_{w_{n}}^{u_{n}}\right) =L_{v}.$
\end{definition}

\strut

For convenience, we will define the following subset of $NC(n)$ ;

\strut

\begin{definition}
Let $NC(n)$ be the set of all noncrossing partitions over $\{1,...,n\}$ and
fix a $D_{G}$-valued random variable $L_{w_{1}}^{u_{1}}...L_{w_{n}}^{u_{n}}$
in $\left( W^{*}(G),E\right) ,$ where $u_{1},$ ..., $u_{n}\in \{1,*\}.$ For
the fixed $D_{G}$-valued random variable $%
L_{w_{1}}^{u_{1}}...L_{w_{n}}^{u_{n}},$define

\strut 

$\ \ \ C_{w_{1},...,w_{n}}^{u_{1},...,u_{n}}\overset{def}{=}\{\pi \in
NC(n):L_{w_{1}}^{u_{1}}...L_{w_{n}}^{u_{n}}$ is $\pi $-connected$\},$

\strut 

in $NC(n).$ Let $\mu $ be the M\"{o}bius function in the incidence algebra $%
I_{2}.$ Define the number $\mu _{w_{1},...,w_{n}}^{u_{1},...,u_{n}},$ for
the fixed $D_{G}$-valued random variable $%
L_{w_{1}}^{u_{1}}...L_{w_{n}}^{u_{n}},$ by

\strut 

$\ \ \ \ \ \ \ \ \ \ \ \ \ \ \ \mu _{w_{1},...,w_{n}}^{u_{1},...,u_{n}}%
\overset{def}{=}\underset{\pi \in C_{w_{1},...,w_{n}}^{u_{1},...,u_{n}}}{%
\sum }\mu (\pi ,1_{n}).$
\end{definition}

\bigskip \strut

Assume that there exists $\pi \in NC(n)$ such that $%
L_{w_{1}}^{u_{1}}...L_{w_{n}}^{u_{n}}=L_{v}$ is $\pi $-connected. Then $\pi
\in C_{w_{1},...,w_{n}}^{u_{1},...,u_{n}}$ and there exists the maximal
partition $\pi _{0}\in C_{w_{1},...,w_{n}}^{u_{1},...,u_{n}}$ such that $%
L_{w_{1}}^{u_{1}}...L_{w_{n}}^{u_{n}}=L_{v}$ is $\pi _{0}$-connected.
(Recall that $NC(n)$ is a lattice. We can restrict this lattice ordering on $%
C_{w_{1},...,w_{n}}^{u_{1},...,u_{n}}$ and hence it is a POset, again.)
Notice that $1_{n}\in C_{w_{1},...,w_{n}}^{u_{1},...,u_{n}}.$ Therefore, the
maximal partition in $C_{w_{1},...,w_{n}}^{u_{1},...,u_{n}}$ is $1_{n}.$
Hence we have that ;

\bigskip \strut

\begin{lemma}
Let $L_{w_{1}}^{u_{1}}...L_{w_{n}}^{u_{n}}\in \left( W^{*}(G),E\right) $ be
a $D_{G}$-valued random variable having the $*$-axis-property. Then

\strut 

$\ \ \ \ \ \ \ \ \ \ \ E\left( L_{w_{1}}^{u_{1}}...L_{w_{n}}^{u_{n}}~\right)
=\ \widehat{E}(\pi )\left( L_{w_{1}}^{u_{1}}~\otimes ...\otimes
~L_{w_{n}}^{u_{n}}\right) ,$

\strut 

for all $\pi \in C_{w_{1},...,w_{n}}^{u_{1},...,u_{n}}.$
\end{lemma}

\bigskip

\begin{proof}
By the previous discussion, we can get the result.
\end{proof}

\strut \strut

By the previous lemmas, we have that

\strut

\begin{theorem}
Let $n\in 2\mathbb{N}$ and let $L_{w_{1}}^{u_{1}},...,L_{w_{n}}^{u_{n}}\in
\left( W^{*}(G),E\right) $ be $D_{G}$-valued random variables, where $%
w_{1},...,w_{n}\in FP(G)$ and $u_{j}\in \{1,*\},$ $j=1,...,n.$ Then

\strut 

$\ \ \ \ \ \ \ k_{n}\left( L_{w_{1}}^{u_{1}}...L_{w_{n}}^{u_{n}}~\right)
=\mu _{w_{1},...,w_{n}}^{u_{1},...,u_{n}}\cdot
E(L_{w_{1}}^{u_{1}},...,L_{w_{n}}^{u_{n}}),$

\strut 

where $\mu _{w_{1},...,w_{n}}^{u_{1},...,u_{n}}=\underset{\pi \in
C_{w_{1},...,w_{n}}^{u_{1},...,u_{n}}}{\sum }\mu (\pi ,1_{n}).$
\end{theorem}

\bigskip \strut \strut

\begin{proof}
We can compute that

\strut

$\ k_{n}\left( L_{w_{1}}^{u_{1}},\,...,L_{w_{n}}^{u_{n}}\right) =\underset{%
\pi \in NC(n)}{\sum }\widehat{E}(\pi )\left( L_{w_{1}}^{u_{1}}~\otimes
...\otimes ~L_{w_{n}}^{u_{n}}\right) \mu (\pi ,1_{n})$

\strut

$\ \ \ \ \ =\underset{\pi \in C_{w_{1},...,w_{n}}^{u_{1},...,u_{n}}}{\sum }%
\widehat{E}(\pi )\left( L_{w_{1}}^{u_{1}}~\otimes ...\otimes
~L_{w_{n}}^{u_{n}}\right) \mu (\pi ,1_{n})$

\strut

by the $\pi $-connectedness

\strut

$\ \ \ \ \ \ =\underset{\pi \in C_{w_{1},...,w_{n}}^{u_{1},...,u_{n}}}{\sum }%
E\left( L_{w_{1}}^{u_{1}}...~L_{w_{n}}^{u_{n}}\right) \mu (\pi ,1_{n})$

\strut

by the previous lemma

\strut

$\ \ \ \ \ \ =\left( \underset{\pi \in C_{w_{1},...,w_{n}}^{u_{1},...,u_{n}}%
}{\sum }\mu (\pi ,1_{n})\right) \cdot E\left(
L_{w_{1}}^{u_{1}}...~L_{w_{n}}^{u_{n}}\right) .$

\strut
\end{proof}

\strut

Now, we can get the following $D_{G}$-valued cumulants of the random
variable $a$ ;

\strut

\begin{corollary}
Let $n\in \mathbb{N}$ and let $a=a_{d}+a_{(*)}+a_{(non-*)}\in \left(
W^{*}(G),E\right) $ be our $D_{G}$-valued random variable. Then $k_{1}\left(
d_{1}a\right) =d_{1}a_{d}$ and $k_{n}\left( d_{1}a,...,d_{n}a\right)
=0_{D_{G}},$ for all odd $n.$ If $n\in \Bbb{N}\,\setminus \,\{1\},$ then

\strut \strut 

$\ \ k_{n}\left( d_{1}a,...,d_{n}a\right) =~\underset{(v_{1},...,v_{n})\in
\Pi _{j=1}^{n}V(G:d_{j})}{\sum }\left( \Pi _{j=1}^{n}q_{v_{j}}\right) $

\strut 

$\ \ \ \ \ \ \ \ \underset{(w_{1},...,w_{n})\in
FP_{*}(G:a)^{n},\,w_{j}=x_{j}w_{j}y_{j},\,u_{w_{j}}\in
\{1,*\},\,l_{w_{1},...,w_{n}}^{u_{1},...,u_{n}}\in LP_{n}^{*}}{\sum }\left(
\Pi _{j=1}^{n}p_{w_{j}}\right) $

\strut 

$\ \ \ \ \ \ \ \ \ \ \left( \Pi _{j=1}^{n}\delta
_{(v_{j},x_{j},y_{j}:u_{j})}\right) \left( \mu
_{w_{1},...,w_{n}}^{u_{1},...,u_{n}}\cdot E\left(
L_{w_{1}}^{u_{w_{1}}}...L_{w_{n}}^{u_{w_{n}}}\right) \right) ,$

\strut 

where $d_{1},...,d_{n}\in D_{G}$ are arbitrary. $\square $
\end{corollary}

\strut \strut \strut \strut

We have the following trivial $D_{G}$-valued moments and cumulants of an
arbitrary $D_{G}$-valued random variable ;

\strut

\begin{corollary}
Let $a\in \left( W^{*}(G),E\right) $ be a $D_{G}$-valued random variable and
let $n\in \mathbb{N}.$ Then

\strut 

(1) The $n$-th trivial $D_{G}$-valued moment of $a$ is

\strut 

$\ \ \ \ \ E(a^{n})=\,\underset{(w_{1},...,w_{n})\in \mathbb{F}%
^{+}(G:a)^{n},\,u_{w_{j}}\in
\{1,*\},\,l_{w_{1},...,w_{n}}^{u_{1},...,u_{n}}\in LP_{n}^{*}}{\sum }$

\strut 

$\ \ \ \ \ \ \ \ \ \ \ \ \ \ \ \ \ \ \ \ \ \ \ \ \ \ \ \ \ \ \ \ \ \ \ \
\left( \Pi _{j=1}^{n}p_{w_{j}}\right) \cdot E\left(
L_{w_{1}}^{u_{1}},...,~L_{w_{n}}^{u_{n}}\right) .$

\strut 

(2) The $n$-th trivial $D_{G}$-valued cumulant of $a$ is

\strut 

$\ \ \ \ k_{1}\left( a\right) =E(a)=a_{d}$

\strut \strut 

and

\strut 

$\ \ \ \ \ k_{n}\left( \underset{n-times}{\underbrace{a,.....,a}}\right) =~%
\underset{(w_{1},...,w_{n})\in FP_{*}(G:a)^{n},\,\,\,u_{w_{j}}\in
\{1,*\},\,l_{w_{1},...,w_{n}}^{u_{1},...,u_{n}}\in LP_{n}^{*}}{\sum }$

\strut \strut 

$\ \ \ \ \ \ \ \ \ \ \ \ \ \ \ \ \left( \Pi _{j=1}^{n}p_{w_{j}}\right)
\left( \mu _{w_{1},...,w_{n}}^{u_{1},...,u_{n}}\cdot E\left(
L_{w_{1}}^{u_{1}},...,~L_{w_{n}}^{u_{n}}\right) \right) ,$

\strut \strut 

where $d_{1},...,d_{n}\in D_{G}$ are arbitrary. \ $\square $
\end{corollary}

\strut \bigskip

\strut

\strut

\section{\strut $D_{G}$-Freeness on $\left( W^{*}(G),E\right) $}

\strut

\strut

Like before, throughout this chapter, let $G$ be a countable directed graph
and $\left( W^{*}(G),E\right) $, the graph $W^{*}$-probability space over
its diagonal subalgebra $D_{G}.$ In this chapter, we will consider the $%
D_{G} $-valued freeness of given two random variables in $\left(
W^{*}(G),E\right) $. We will characterize the $D_{G}$-freeness of $D_{G}$%
-valued random variables $L_{w_{1}}$ and $L_{w_{2}},$ where $w_{1}\neq
w_{2}\in FP(G).$ And then we will observe the $D_{G}$-freeness of arbitrary
two $D_{G}$-valued random variables $a_{1}$ and $a_{2}$ in terms of their
supports. Let

\strut

(3.0) $\ a=\underset{w\in \mathbb{F}^{+}(G:a),\,u_{w}\in \{1,*\}}{\sum }%
p_{w}L_{w}^{u_{w}}$ \& $b=\underset{w^{\prime }\in \mathbb{F}%
^{+}(G:b),\,u_{w^{\prime }}\in \{1,*\}}{\sum }p_{w^{\prime }}L_{w^{\prime
}}^{u_{w^{\prime }}}$

\strut

be fixed $D_{G}$-valued random variables in $\left( W^{*}(G),E\right) $.

\strut

Now, fix $n\in \mathbb{N}$ and let $\left( a_{i_{1}}^{\varepsilon
_{i_{1}}},...,a_{i_{n}}^{\varepsilon _{i_{n}}}\right) \in
\{a,b,a^{*},b^{*}\}^{n},$ where $\varepsilon _{i_{j}}\in \{1,*\}.$ For
convenience, put

\strut

\begin{center}
$a_{i_{j}}^{\varepsilon _{i_{j}}}=\underset{w_{i_{j}}\in \mathbb{F}%
^{+}(G:a),\,\,u_{j}\in \{1,\ast \}}{\sum }p_{w_{j}}^{(j)}L_{w_{j}}^{u_{j}},$
for \ $j=1,...,n.$
\end{center}

\strut

Then, by the little modification of Section , we have that ;

\strut

(3.1)

\strut

$\ E\left( d_{i_{1}}a_{i_{1}}^{\varepsilon
_{i_{1}}}...d_{i_{n}}a_{i_{n}}^{\varepsilon _{i_{n}}}\right) $

\strut

$\ \ \ \ \ =\,\underset{(v_{i_{1}},...,v_{i_{n}})\in \Pi
_{k=1}^{n}V(G:d_{i_{k}})}{\sum }\left( \Pi _{k=1}^{n}q_{v_{i_{k}}}\right) $

\strut

$\ \ \ \ \ \ \ \ \ \ \ \ \ \ \ \ \ \underset{(w_{i_{1}},...,w_{i_{n}})\in
\Pi _{k=1}^{n}\mathbb{F}^{+}(G:a_{i_{k}}),%
\,w_{i_{j}}=x_{i_{j}}w_{i_{j}}y_{i_{j}},u_{i_{j}}\in \{1,\ast \}}{\sum }%
\left( \Pi _{k=1}^{n}p_{w_{i_{k}}}^{(k)}\right) $

\strut

$\ \ \ \ \ \ \ \ \ \ \ \ \ \ \ \ \ \ \ \ \ \ \ \ \ \ \ \ \ \ \ \ \ \ \left(
\Pi _{j=1}^{n}\delta _{(v_{i_{j}},x_{i_{j}},y_{i_{j}}:u_{i_{j}})}\right)
E\left( L_{w_{i_{1}}}^{u_{i_{1}}}...L_{w_{i_{n}}}^{u_{i_{n}}}\right) .$

\strut

\strut

Therefore, we have that

\strut

(3.2)

\strut

$\ \ k_{n}\left( d_{i_{1}}a_{i_{1}}^{\varepsilon
_{i_{1}}},...,d_{i_{n}}a_{i_{n}}^{\varepsilon _{i_{n}}}\right) =\underset{%
(v_{i_{1}},...,v_{i_{n}})\in \Pi _{k=1}^{n}V(G:d_{i_{k}})}{\sum }\left( \Pi
_{k=1}^{n}q_{v_{i_{k}}}\right) $

\strut

$\ \ \ \ \ \ \ \ \ \ \ \ \ \ \ \ \ \underset{(w_{i_{1}},...,w_{i_{n}})\in
\Pi
_{k=1}^{n}FP(G:a_{i_{k}}),\,w_{i_{j}}=x_{i_{j}}w_{i_{j}}y_{i_{j}},u_{i_{j}}%
\in \{1,\ast \}}{\sum }\left( \Pi _{k=1}^{n}p_{w_{i_{k}}}^{(k)}\right) $

\strut

$\ \ \ \ \ \ \ \ \ \ \ \ \ \ \ \ \ \ \ \ \left( \Pi _{j=1}^{n}\delta
_{(v_{i_{j}},x_{i_{j}},y_{i_{j}}:u_{i_{j}})}\right) \left( \mu
_{w_{i_{1}},...,w_{i_{n}}}^{u_{i_{1}},...,u_{i_{n}}}\cdot E\left(
L_{w_{i_{1}}}^{u_{i_{1}}}...L_{w_{i_{n}}}^{u_{i_{n}}}\right) \right) ,$

\strut

where $\mu _{w_{1},...,w_{n}}^{u_{1},...,u_{n}}=\underset{\pi \in
C_{w_{i_{1}},...,w_{i_{n}}}^{u_{i_{1}},...,u_{i_{n}}}}{\sum }\mu (\pi
,1_{n}) $ and

\bigskip

\begin{center}
$C_{w_{i_{1}},...,w_{i_{n}}}^{u_{i_{1}},...,u_{i_{n}}}=\{\pi \in
NC^{(even)}(n):L_{w_{1}}^{u_{1}}...L_{w_{n}}^{u_{n}}$ is $\pi $-connected$%
\}. $
\end{center}

\strut

So, we have the following proposition, by the straightforward computation ;

\bigskip \strut \strut

\begin{proposition}
Let $a,b\in \left( W^{*}(G),E\right) $ be $D_{G}$-valued random variables,
such that $a\notin W^{*}(\{b\},D_{G}),$ and let $\left(
a_{i_{1}}^{\varepsilon _{i_{1}}},...,a_{i_{n}}^{\varepsilon _{i_{n}}}\right)
\in \{a,b,a^{*},b^{*}\}^{n},$ for $n\in \mathbb{N}\setminus \{1\},$ where $%
\varepsilon _{i_{j}}\in \{1,*\},$ $j=1,...,n.$ Then\strut 

\strut 

\strut (3.3)

\strut 

$\ \ \ \ \ k_{n}\left( d_{i_{1}}a_{i_{1}}^{\varepsilon
_{i_{1}}},...,d_{i_{n}}a_{i_{n}}^{\varepsilon _{i_{n}}}\right) $

\strut 

$\ \ \ \ \ \ =\underset{(v_{1},...,v_{n})=(x,y,...,x,y)\in \Pi
_{j=1}^{n}V(G:d_{j})}{\sum }\left( \Pi _{j=1}^{n}q_{v_{j}}\right) $

\strut 

$\ \ \ \ \ \ \ \ \ \ \underset{(w_{i_{1}},...,w_{i_{n}})\in \left( \Pi
_{k=1}^{n}FP_{*}(G:a_{i_{k}})\right) \cup
W_{*}^{iw_{1},...,w_{n}},\,w_{i_{j}}=x_{i_{j}}w_{i_{j}}y_{i_{j}}}{\sum }%
\left( \Pi _{k=1}^{n}p_{w_{i_{j}}}^{(k)}\right) $

\strut 

$\ \ \ \ \ \ \ \ \ \ \ \ \ \ \left( \Pi _{j=1}^{n}\delta
_{(v_{i_{j}},x_{i_{j}},y_{i_{j}}:u_{i_{j}})}\right) \left( \mu
_{w_{i_{1}},...,w_{i_{n}}}^{u_{i_{1}},...,u_{i_{n}}}\cdot \Pr oj\left(
L_{w_{i_{1}}}^{u_{i_{1}}}...L_{w_{i_{n}}}^{u_{i_{n}}}\right) \right) $

\strut 

\strut where $\mu _{n}=\underset{\pi \in
C_{w_{i_{1}},...,w_{i_{n}}}^{u_{i_{1}},...,u_{i_{n}}}}{\sum }\mu (\pi ,1_{n})
$ and

\strut 

$\ \ \ \ \ \ W_{*}^{w_{1},...,w_{n}}=\{w\in FP_{*}^{c}(G:a)\cup
FP_{*}^{c}(G:b):$

$\ \ \ \ \ \ \ \ \ \ \ \ \ \ \ \ \ \ \ \ \ \ \ \ \ \ \ \ \ \ \ \ \ $both $%
L_{w}^{u_{w}}$ and $L_{w}^{u_{w}\,*}$\ are in $%
L_{w_{1}}^{u_{w_{1}}}...L_{w_{n}}^{u_{w_{n}}}\}.$

$\square $
\end{proposition}

\strut \strut \strut \strut \strut \strut

\begin{corollary}
Let $x$ and $y$ be the $D_{G}$-valued random variables in $\left(
W^{*}(G),E\right) $. The $D_{G}$-valued random variables $a$ and $b$ are
free over $D_{G}$ in $\left( W^{*}(G),E\right) $ if

\strut \strut \strut 

$\ \ \ \ \ \ \ FP_{*}\left( G:P(x,x^{*})\right) \cap FP_{*}\left(
G:Q(y,y^{*})\right) =\emptyset $

and

$\ \ \ \ \ \ \ \ \ \ \ \ \ \ \ W_{*}^{\{P(x,x^{*}),\,Q(y,y^{*})\}}=\emptyset
,$

\strut 

for all $P,Q\in \mathbb{C}[z_{1},z_{2}].$ \ $\square $
\end{corollary}

\strut \strut \strut \strut

By using (3.2), we can compute the mixed $D_{G}$-valued cumulants of two $%
D_{G}$-valued random variables. However, the formula is very abstract. So,
we will consider the above formula for fixed two generators of $W^{\ast }(G)$%
.

\strut

\begin{definition}
Let $G$ be a countable directed graph and $\mathbb{F}^{+}(G)$, the free
semigroupoid of $G$ and let $FP(G)$ be the subset of $\mathbb{F}^{+}(G)$
consisting of all finite paths. Define a subset $loop(G)$ of $FP(G)$
containing all loop finite paths or loops. (Remark that, in general, loop
finite paths are different from loop-edges. Clearly, all loop-edges are
loops in $FP(G).$) i.e.,

\strut 

$\ \ \ \ \ \ \ loop(G)\overset{def}{=}\{l\in FP(G):l$ is a loop$\}\subset
FP(G).$

\strut 

Also define the subset $loop^{c}(G)$ of $FP(G)$ consisting of all non-loop
finite path by

\strut 

$\ \ \ \ \ \ \ \ \ loop^{c}(G)\overset{def}{=}FP(G)\,\setminus \,loop(G).$

\strut 

Let $l\in loop(G)$ be a loop finite path. We say that $l$ is a \textbf{basic
loop} if there exists no loop $w\in loop(G)$ such that $l=w^{k},$ $k\in %
\mathbb{N}\,\setminus \,\{1\}.$ Define

\strut 

$\ \ \ Loop(G)\overset{def}{=}\{l\in loop(G):l$ is a basic loop$%
\}\subsetneqq loop(G).$

\strut 

Let $l_{1}=w_{1}^{k_{1}}$ and $l_{2}=w_{2}^{k_{2}}$ in $loop(G),$ where $%
w_{1},w_{2}\in Loop(G).$ We will say that the loops $l_{1}$ and $l_{2}$ are 
\textbf{diagram-distinct} if $w_{1}\neq w_{2}$ in $Loop(G).$ Otherwise, they
are not diagram-distinct.
\end{definition}

\strut

Now, we will introduce the more general diagram-distinctness of general
finite paths ;

\strut

\begin{definition}
(\textbf{Diagram-Distinctness}) We will say that the finite paths $w_{1}$
and $w_{2}$ are \textbf{diagram-distinct} if $w_{1}$ and $w_{2}$ have
different diagrams in the graph $G.$ Let $X_{1}$ and $X_{2}$ be subsets of $%
FP(G).$ The subsets $X_{1}$ and $X_{2}$ are said to be diagram-distinct if $%
x_{1}$ and $x_{2}$ are diagram-distinct for all pairs $(x_{1},x_{2})$ $\in $ 
$X_{1}\times X_{2}.$ This diagram-distinctness implies the
diagram-distinctness of loops.
\end{definition}

\strut \strut

Let $H$ be a directed graph with $V(H)=\{v_{1},v_{2}\}$ and $%
E(H)=\{e_{1}=v_{1}e_{1}v_{2},e_{2}=v_{2}e_{2}v_{1}\}.$ Then $l=e_{1}e_{2}$
is a loop in $FP(H)$ (i.e., $l\in loop(H)$). Moreover, it is a basic loop
(i.e., $l\in Loop(H)$). However, if we have a loop $%
w=e_{1}e_{2}e_{1}e_{2}=l^{2},$ then it is not a basic loop. i.e.,

\strut

\begin{center}
$l^{2}\in loop(H)\,\setminus \,Loop(H).$
\end{center}

\strut

If the graph $G$ contains at least one basic loop $l\in FP(G)$, then we have

\strut

\begin{center}
$\{l^{n}:n\in \mathbb{N}\}\subseteq loop(G)$ \ and \ $\{l\}\subseteq
Loop(G). $
\end{center}

\strut

Suppose that $l_{1}$ and $l_{2}$ are not diagram-distinct. Then, by
definition, there exists $w\in Loop(G)$ such that $l_{1}=w^{k_{1}}$ and $%
l_{2}=w^{k_{2}},$ for some $k_{1},k_{2}\in \mathbb{N}.$ On the graph $G,$
indeed, $l_{1}$ and $l_{2}$ make the same diagram. On the other hands, we
can see that if $w_{1}\neq w_{2}\in loop^{c}(G),$ then they are
automatically diagram-distinct.

\strut

\begin{lemma}
Suppose that $w_{1}\neq w_{2}\in loop^{c}(G)$ with $w_{1}=v_{11}w_{1}v_{12}$
and $w_{2}=v_{21}w_{2}v_{22}.$ Then $L_{w_{1}}$ and $L_{w_{2}}$ are free
over $D_{G}$ in $\left( W^{*}(G),E\right) .$
\end{lemma}

\strut

\begin{proof}
By definition, $L_{w_{1}}$ and $L_{w_{2}}$ are free over $D_{G}$ if and only
if all mixed $D_{G}$-valued cumulants of $W^{\ast }(\{L_{w_{1}}\},D_{G})$
and $W^{\ast }(\{L_{w_{2}}\},D_{G})$ vanish. Equivalently, all $D_{G}$%
-valued cumulants of $P\left( L_{w_{1}},L_{w_{1}}^{\ast }\right) $ and $%
Q\left( L_{w_{2}},L_{w_{2}}^{\ast }\right) $ vanish, for all $P,Q\in \mathbb{%
C}[z_{1},z_{2}].$ Since $w_{1}\neq w_{2}$ are non-loop edges, we can easily
verify that $w_{1}^{k_{1}}$ and $w_{2}^{k_{2}}$ are not admissible (i.e., $%
w_{1}^{k_{1}}\notin \mathbb{F}^{+}(G)$ and $w_{1}^{k_{2}}\notin \mathbb{F}%
^{+}(G)$), for all $k_{1},k_{2}\in \mathbb{N}\,\setminus \,\{1\}.$ This
shows that

\strut

$\ \ \ \ \ \ \ \ \ L_{w_{j}}^{k}=0_{D_{G}}=\left( L_{w_{j}}^{k\,}\right)
^{*},$ for $j=1,2.$

\strut

Thus, to show that $L_{w_{1}}$ and $L_{w_{2}}$ are free over $D_{G},$ it
suffices to show that all mixed $D_{G}$-valued cumulants of $P\left(
L_{w_{1}},L_{w_{1}}^{*}\right) $ and $Q\left( L_{w_{2}},L_{w_{2}}^{*}\right) 
$ vanish, for all $P,Q\in \mathbb{C}[z_{1},z_{2}]$ such that

\strut

$\ \ \ \ \ \ \ \ P(z_{1},z_{2})=\alpha _{1}z_{1}+\alpha
_{2}z_{1}z_{2}+\alpha _{3}z_{2}z_{1}+\alpha _{4}z_{2}$

and

$\ \ \ \ \ \ \ \ Q(z_{1},z_{2})=\beta _{1}z_{1}+\beta _{2}z_{1}z_{2}+\beta
_{3}z_{2}z_{1}+\beta _{4}z_{2},$

\strut \strut

where $\alpha ,\beta \in \mathbb{C}.$ So, for such $P$ and $Q,$ we have that

\strut

\ \ \ \ \ \ \strut $P\left( L_{w_{1}},L_{w_{1}}^{*}\right) =\alpha
_{1}L_{w_{1}}+\alpha _{2}L_{v_{11}}+\alpha _{3}L_{v_{12}}+\alpha
_{4}L_{w_{1}}^{*}$

and

$\ \ \ \ \ \ Q\left( L_{w_{2}},L_{w_{2}}^{*}\right) =\beta
_{1}L_{w_{2}}+\beta _{2}L_{v_{21}}+\beta _{3}L_{v_{22}}+\beta
_{4}L_{w_{2}}^{*}.$

\strut

Thus, we have that

\strut

$\ \ FP_{*}\left( G:P(L_{w_{1}},L_{w_{1}}^{*})\right) \supseteq \{w_{1}\},$
\ $FP_{*}\left( G:Q(L_{w_{2}},L_{w_{2}}^{*})\right) \supseteq \{w_{2}\}$

\strut

and

\strut

$\ \ FP_{*}^{c}\left( G:P(L_{w_{1}},L_{w_{1}}^{*})\right) \supseteq
\{w_{1}\},$ \ $FP_{*}^{c}\left( G:Q(L_{w_{2}},L_{w_{2}}^{*})\right)
\supseteq \{w_{2}\}.$

\strut

Remark that if $FP_{*}\left( G:P(L_{w_{1}},L_{w_{1}}^{*})\right) =\{w_{1}\},$
then $FP_{*}^{c}\left( G:P(L_{w_{1}},L_{w_{1}}^{*})\right) =\emptyset $, and
if$\ FP_{*}^{c}\left( G:P(L_{w_{1}},L_{w_{1}}^{*})\right) =\{w_{1}\},$ then $%
FP_{*}\left( G:P(L_{w_{1}},L_{w_{1}}^{*})\right) =\emptyset .$ The similar
relation holds for $Q\left( L_{w_{2}},L_{w_{2}}^{*}\right) .$ So, we have
that

\strut

$\ \ \ \ \ \ FP_{*}\left( G:P(L_{w_{1}},L_{w_{1}}^{*})\right) \cap
FP_{*}\left( G:Q(L_{w_{2}},L_{w_{2}}^{*})\right) =\emptyset $

and

$\ \ \ \ \ \ \ \ \ \ \ \ \
W_{*}^{\{P(L_{w_{1}},L_{w_{1}}^{*}),\,Q(L_{w_{2}},L_{w_{2}}^{*})\}}=%
\emptyset .$

\strut

Therefore, by the formula (3.4.3), we have the vanishing mixed $D_{G}$%
-valued cumulants of $P\left( L_{w_{1}},L_{w_{1}}^{\ast }\right) $ and $%
Q\left( L_{w_{2}},L_{w_{2}}^{\ast }\right) ,$ for all $n\in \mathbb{N}$ and
for all such $P,Q\in \mathbb{C}[z_{1},z_{2}].$ So, we can conclude that $%
L_{w_{1}}$ and $L_{w_{2}}$ are free over $D_{G}$ in $\left( W^{\ast
}(G),E\right) .$\strut
\end{proof}

\strut \strut \strut \strut \strut \strut \strut \strut \strut \strut

Now, we will consider the loop case.

\strut \strut

\begin{lemma}
Let $l_{1}\neq l_{2}\in Loop(G)$ be\textbf{\ basic loops} such that $%
l_{1}=v_{1}l_{1}v_{1}$ and $l_{2}=v_{2}l_{2}v_{2},$ for $v_{1},v_{2}\in V(G)$
(possibly $v_{1}=v_{2}$). i.e.e, two basic loops $l_{1}$ and $l_{2}$ are
diagram-distinct. Then the $D_{G}$-valued random variables $L_{l_{1}}$ and $%
L_{l_{2}}$ are free over $D_{G}$ in $\left( W^{*}(G),E\right) .$
\end{lemma}

\strut

\begin{proof}
Different from the non-loop case, if $l_{1}$ and $l_{2}$ are loops, then $%
l_{1}^{k_{1}}$ and $l_{2}^{k_{2}}$ exist in $FP(G),$ for all $k_{1},k_{2}\in %
\mathbb{N}.$ To show that $L_{l_{1}}$ and $L_{l_{2}}$ are free over $D_{G},$
it suffices to show that all mixed $D_{G}$-valued cumulants of $P\left(
L_{w_{1}},L_{w_{1}}^{\ast }\right) $ and $Q\left( L_{w_{2}},L_{w_{2}}^{\ast
}\right) $ vanish, for all $P,Q\in \mathbb{C}[z_{1},z_{2}].$ such that

\strut

$\ \ \ \ \ \ \ \ P(z_{1},z_{2})=f_{1}(z_{1})+f_{2}(z_{2})+P_{0}(z_{1},z_{2})$

and

$\ \ \ \ \ \ \ \
Q(z_{1},z_{2})=g_{1}(z_{1})+g_{2}(z_{2})+Q_{0}(z_{1},z_{2}), $

\strut \strut

where $f_{1},f_{2},g_{1},g_{2}\in \mathbb{C}[z]$ and $P_{0},Q_{0}\in \mathbb{%
C}[z_{1},z_{2}]$ such that $P_{0}$ and $Q_{0}$ does not contain polynomials
only in $z_{1}$ and $z_{2}.$ So, for such $P$ and $Q,$ we have that

\strut

\ \ \ \ \ \ \ \strut $P\left( L_{l_{1}},L_{l_{1}}^{*}\right)
=f_{1}(L_{l_{1}})+f_{2}(L_{l_{1}}^{*})+P_{0}(L_{l_{1}},L_{l_{1}}^{*})$

and

$\ \ \ \ \ \ \ Q\left( L_{l_{2}},L_{l_{2}}^{*}\right)
=g_{1}(L_{l_{2}})+g_{2}(L_{l_{2}}^{*})+Q_{0}(L_{l_{2}},L_{l_{2}}^{*}).$

\strut

Notice that $L_{l_{j}}^{k}=L_{l_{j}^{k}},$ for all $k\in \mathbb{N},$ \ $%
j=1,2.$ Also, notice that

\strut

$\ \ \ \ \ \ \ \ \ P_{0}\left( L_{l_{1}},L_{l_{1}}^{*}\right)
=f_{1}^{0}(L_{l_{1}})+f_{2}^{0}(L_{l_{1}}^{*})+\alpha L_{v_{1}}$

and

$\ \ \ \ \ \ \ \ \ Q_{0}\left( L_{l_{2}},L_{l_{2}}^{*}\right)
=g_{1}^{0}(L_{l_{2}})+g_{2}^{0}(L_{l_{2}}^{*})+\beta L_{v_{2}},$

\strut

where $f_{1}^{0},f_{2}^{0},g_{1}^{0},g_{2}^{0}\in \mathbb{C}[z]$ and $\alpha
,\beta \in \mathbb{C}$, by the fact that

\strut

\ \ \ \ \ \ \ \ \ \ \ \ \ \ \ \ \ $\ \
L_{l_{j}}^{*}L_{l_{j}}=L_{v_{j}}=L_{l_{j}}L_{l_{j}}^{*},$

\strut

under the weak-topology. So, finally, we have that

\strut

$\ \ \ \ \ \ \ \ \ P\left( L_{l_{1}},L_{l_{1}}^{*}\right) =\mathbf{f}%
_{1}(L_{l_{1}})+\mathbf{f}_{2}(L_{l_{1}}^{*})+\alpha L_{v_{1}}$

and

$\ \ \ \ \ \ \ \ \ Q\left( L_{l_{2}},L_{l_{2}}^{*}\right) =\mathbf{g}%
_{1}(L_{l_{2}})+\mathbf{g}_{1}(L_{l_{2}}^{*})+\beta L_{v_{2}},$

\strut

where $\mathbf{f}_{1},\mathbf{f}_{2},\mathbf{g}_{1},\mathbf{g}_{2}\in %
\mathbb{C}[z]$ and $\alpha ,\beta \in \mathbb{C}.$ Thus, we have that

\strut

$\ \ FP_{*}\left( G:P(L_{l_{1}},L_{l_{1}}^{*})\right) \subseteq
\{l_{1}^{k}\}_{k=1}^{\infty },$ \ $FP_{*}\left(
G:Q(L_{l_{2}},L_{l_{2}}^{*})\right) \subseteq \{l_{2}^{k}\}_{k=1}^{\infty }$

\strut

and

\strut

$\ \ FP_{*}^{c}\left( G:P(L_{l_{1}},L_{l_{1}}^{*})\right) \subseteq
\{l_{1}^{k}\}_{k=1}^{\infty },$ \ $FP_{*}^{c}\left(
G:Q(L_{l_{2}},L_{l_{2}}^{*})\right) \subseteq \{l_{2}^{k}\}_{k=1}^{\infty }.$

\strut

So, we have that

\strut

$\ \ \ \ \ \ FP_{*}\left( G:P(L_{w_{1}},L_{w_{1}}^{*})\right) \cap
FP_{*}\left( G:Q(L_{w_{2}},L_{w_{2}}^{*})\right) =\emptyset ,$

\strut

because $l_{1}$ and $l_{2}$ are in $Loop(G)$ (and hence if $l_{1}\neq l_{2},$
then they are diagram-distinct.) And we have that

\strut

$\ \ \ \ \ \ \ \ \ \ \ \ \
W_{*}^{\{P(L_{w_{1}},L_{w_{1}}^{*}),\,Q(L_{w_{2}},L_{w_{2}}^{*})\}}=%
\emptyset .$

\strut

Therefore, by the formula (3.4.3), we have the vanishing mixed $D_{G}$%
-valued cumulants of $P\left( L_{l_{1}},L_{l_{1}}^{\ast }\right) $ and $%
Q\left( L_{l_{2}},L_{l_{2}}^{\ast }\right) ,$ for all $n\in \mathbb{N}$ and
for all $P,Q\in \mathbb{C}[z_{1},z_{2}].$ Since $P$ and $Q$ are arbitrary,
we can conclude that $L_{l_{1}}$ and $L_{l_{2}}$ are free over $D_{G}$ in $%
\left( W^{\ast }(G),E\right) .$\strut
\end{proof}

\strut

Notice that we assumed that the loops $l_{1}$ and $l_{2}$ are basic loops in
the previous lemma. Since they are distinct basic loops, they are
automatically diagram-distinct. Now, assume that $l_{1}$ and $l_{2}$ are not
diagram-distinct. i.e., there exists a basic loop $w\in Loop(G)$ such that $%
l_{1}=w^{k_{1}}$ and $l_{2}=w^{k_{2}},$ for some $k_{1},k_{2}\in \mathbb{N}.$
In other words, the loops $l_{1}$ and $l_{2}$ have the same diagram in the
graph $G.$ Then the $D_{G}$-valued random variables $L_{l_{1}}$ and $%
L_{l_{2}}$ are not free over $D_{G}$ in $\left( W^{\ast }(G),E\right) $. See
the next example ;

\strut

\begin{example}
Let $G_{1}$ be a directed graph with $V(G_{1})=\{v\}$ and $%
E(G_{1})=\{l=vlv\}.$ So, in this case,

\strut 

$\ \ \ \ \ \ \ E(G_{1})=Loop(G_{1}),$ \ $FP(G_{1})=loop(G_{1}),$

\strut and

$\ \ \ \ \ \ \ \ \ \ \ loop(G_{1})=\{l^{k}:k\in \mathbb{N}\}.$

\strut 

Thus, even if $w_{1}\neq w_{2}\in loop(G_{1}),$ $w_{1}$ and $w_{2}$ are Not
diagram-distinct. Take $l^{2}$ and $l^{3}$ in $FP(G_{1}).$ Then the $%
D_{G_{1}}$-valued random variable $L_{l^{2}}$ and $L_{l^{3}}$ are not free
over $D_{G_{1}}$ in $\left( W^{*}(G_{1}),E\right) .$ Indeed, let's take $%
P,Q\in \mathbb{C}[z_{1},z_{2}]$ as

\strut 

$\ \ \ \ \ \ \ P(z_{1},z_{2})=z_{1}^{3}+z_{2}^{3}$ \ \ and \ \ $%
Q(z_{1},z_{2})=z_{1}^{2}+z_{2}^{2}.$

\strut Then

$\ \ \ \ \ \ \ \ \ P\left( L_{l^{2}},L_{l^{2}}^{*}\right)
=L_{l^{2}}^{3}+L_{l^{2}}^{*\,3}=L_{l^{6}}+L_{l^{6}}^{*}$

and

\ \ \ \ \ \ \ \ \ $Q\left( L_{l^{3}},L_{l^{3}}^{*}\right)
=L_{l^{3}}^{2}+L_{l^{3}}^{*\,\,2}=L_{l^{6}}+L_{l^{6}}^{*}.$

\strut 

Then

\strut 

$\ k_{2}\left(
P(L_{l^{2}},L_{l^{2}}^{*}),\,Q(L_{l^{3}},L_{l^{3}}^{*})\right) =k_{2}\left(
L_{l^{6}}+L_{l^{6}}^{*},\,L_{l^{6}}+L_{l^{6}}^{*}\right) $

\strut 

$\ \ \ =\mu _{l^{6},l^{6}}^{1,*}\Pr oj\left( L_{l^{6}},L_{l^{6}}^{*}\right)
+\mu _{l^{6},l^{6}}^{*,1}\Pr oj\left( L_{l^{6}}^{*},L_{l^{6}}\right) $

\strut 

$\ \ \ =\mu _{l^{6},l^{6}}^{1,*}L_{v}+\mu _{l^{6},l^{6}}^{*,1}L_{v}=\left(
\mu _{l^{6},l^{6}}^{1,*}+\mu _{l^{6},l^{6}}^{*,1}\right) L_{v}$

\strut 

$\ \ \ =2L_{v}\neq 0_{D_{G}},$

\strut 

since $\mu _{l^{6},l^{6}}^{1,*}=\mu (1_{2},1_{2})=1=\mu _{l^{6},l^{6}}^{*,1}.
$ This says that there exists at least one nonvanishing mixed $D_{G}$-valued
cumulant of $W^{*}\left( \{L_{l^{3}}\},D_{G_{1}}\right) $ and $W^{*}\left(
\{L_{l^{2}}\},D_{G_{1}}\right) .$ Therefore, $L_{l^{3}}$ and $L_{l^{2}}$ are
not free over $D_{G_{1}}$ in $\left( W^{*}(G_{1}),E\right) .$
\end{example}

\strut

As we have seen before, if two loops $l_{1}$ and $l_{2}$ are not
diagram-distinct, then $D_{G}$-valued random variables $L_{l_{1}}$ and $%
L_{l_{2}}$ are Not free over $D_{G}.$ However, if $l_{1}$ and $l_{2}$ are
diagram-distinct, we can have the following lemma, by the previous lemma ;

\strut

\begin{lemma}
Let $l_{1}\neq l_{2}\in loop(G)$ be loops and assume that $%
l_{1}=w_{1}^{k_{1}}$ and $l_{2}=w_{2}^{k_{2}},$ where $w_{1},w_{2}\in Loop(G)
$ are basic loops and $k_{1},k_{2}\in \mathbb{N}.$ If $w_{1}\neq w_{2}\in
Loop(G),$ then the $D_{G}$-valued random variables $L_{l_{1}}$ and $L_{l_{2}}
$ are free over $D_{G}$ in $\left( W^{*}(G),E\right) .$ $\square $
\end{lemma}

\strut \strut \strut \strut

Finally, we will observe the following case when we have a loop and a
non-loop finite path.

\strut

\begin{lemma}
Let $l\in loop(G)$ and $w\in loop^{c}(G).$ Then the $D_{G}$-valued random
variables $L_{l}$ and $L_{w}$ are free over $D_{G}$ in $\left(
W^{*}(G),E\right) .$
\end{lemma}

\strut

\begin{proof}
let $l\in loop(G)$ and $w\in loop^{c}(G)$ and let $L_{l}$ and $L_{w}$ be the
corresponding $D_{G}$-valued random variables in $\left( W^{\ast
}(G),E\right) .$ Then, for all $P,Q\in \mathbb{C}[z_{1},z_{2}],$ we have that

\strut

$\ \ \ \ \ \ \ \ FP_{*}\left( G:P(L_{l},L_{l}^{*})\right) \cap FP_{*}\left(
G:Q(L_{w},L_{w}^{*})\right) =\emptyset ,$

\strut

\strut since

$\ \ \ \ \ \ \ \ \ FP_{*}\left( G:P(L_{l},L_{l}^{*})\right) \subseteq
\{l^{k}:k\in \mathbb{N}\}\subset loop(G)$

and

$\ \ \ \ \ \ \ \ \ \ \ \ \ FP_{*}\left( G:Q(L_{w},L_{w}^{*})\right)
=\{w\}\subset loop^{c}(G).$

\strut

Also, since $loop(G)\cap loop^{c}(G)=\emptyset ,$ we can get that

\strut

$\ \ \ \ \ \ \ \ \ \ \ \ \ \ \ \ \
W_{*}^{\{P(L_{l},L_{l}^{*}),\,Q(L_{w},L_{w}^{*})\}}=\emptyset ,$

\strut

for all $P,Q\in \mathbb{C}[z_{1},z_{2}].$ Therefore, the $D_{G}$-valued
random variables $L_{l}$ and $L_{w}$ are free over $D_{G}$ in $\left(
W^{*}(G),E\right) .$
\end{proof}

\strut

\strut Now, we can summarize the above lemmas in this section as follows and
this theorem is one of the main result of this paper. The theorem is the
characterization of $D_{G}$-freeness of generators of $W^{*}(G)$ over $D_{G}$%
.

\strut \strut

\begin{theorem}
Let $w_{1},w_{2}\in FP(G)$ be finite paths. The $D_{G}$-valued random
variables $L_{w_{1}}$ and $L_{w_{2}}$ in $\left( W^{*}(G),E\right) $ are
free over $D_{G}$ if and only if $w_{1}$ and $w_{2}$ are diagram-distinct.
\end{theorem}

\strut

\begin{proof}
$\Leftarrow $) Suppose that finite paths $w_{1}$ and $w_{2}$ are
diagram-distinct. Then the $D_{G}$-valued random variables $L_{w_{1}}$ and $%
L_{w_{2}}$ are free over $D_{G},$ by the previous lemmas.

\strut

($\Rightarrow $) Let $L_{w_{1}}$ and $L_{w_{2}}$ be free over $D_{G}$ in $%
\left( W^{*}(G),E\right) .$ Now, assume that $w_{1}$ and $w_{2}$ are not
diagram-distinct. We will observe the following cases ;

\strut

(Case I) The finite paths $w_{1},w_{2}\in loop(G).$ Since they are not
diagram-distinct, there exists a basic loop $l$ $=$ $vlv$ $\in $ $Loop(G),$
with $v\in V(G),$ such that $w_{1}=l^{k_{1}}$ and $w_{2}=l^{k_{2}},$ for
some $k_{1},k_{2}\in \mathbb{N}.$ As we have seen before, $L_{w_{1}}$ and $%
L_{w_{2}}$ are not free over $D_{G}$ in $\left( W^{*}(G),E\right) .$ Indeed,
if we let $k\in \mathbb{N}$ such that $k_{1}\mid k$ \ and \ $k_{2}\mid k$ \
with $k=k_{1}n_{1}=k_{2}n_{2},$ for $n_{1},n_{2}\in \mathbb{N},$ then we can
take $P,Q\in \mathbb{C}[z_{1},z_{2}]$ defined by

\strut

$\ \ \ \ \ \ P(z_{1},z_{2})=z_{1}^{n_{1}}+z_{2}^{n_{1}}$ \ and \ $%
Q(z_{1},z_{2})=z_{1}^{n_{2}}+z_{2}^{n_{2}}.$

\strut

And then

\strut

$\ \ \ \ P\left( L_{w_{1}},L_{w_{1}}^{*}\right)
=L_{w_{1}}^{n_{1}}+L_{w_{1}}^{*\,%
\,n_{1}}=L_{l^{k_{1}}}^{n_{1}}+L_{l^{k_{1}}}^{*\,%
\,n_{1}}=L_{l^{k}}+L_{l^{k}}^{*}$

\strut

and

\strut

$\ \ \ \ Q\left( L_{w_{2}},L_{w_{2}}^{*}\right)
=L_{w_{2}}^{n_{2}}+L_{w_{2}}^{*\,%
\,n_{2}}=L_{l^{k_{2}}}^{n_{2}}+L_{l^{k_{2}}}^{*%
\,n_{2}}=L_{l^{k}}+L_{l^{k}}^{*}.$

\strut

So,

\strut

$\ \ \ 
\begin{array}{ll}
k_{2}\left( P(L_{w_{1}},L_{w_{1}}^{*}),\,Q(L_{w_{2}},L_{w_{2}}^{*})\right) & 
=k_{2}\left( L_{l^{k}}+L_{l^{k}}^{*},\,L_{l^{k}}+L_{l^{k}}^{*}\right) \\ 
&  \\ 
& =\mu _{l^{k},l^{k}}^{1,*}L_{v}+\mu _{l^{k},l^{k}}^{*,1}L_{v} \\ 
&  \\ 
& =2L_{v}\neq 0_{D_{G}}.
\end{array}
$

\strut

Therefore, $P\left( L_{w_{1}},L_{w_{1}}^{*}\right) $ and $Q\left(
L_{w_{2}},L_{w_{2}}^{*}\right) $ are not free over $D_{G}.$ This shows that $%
W^{*}\left( \{L_{w_{1}}\},D_{G}\right) $ and $W^{*}\left(
\{L_{w_{2}}\},D_{G}\right) $ are not free over $D_{G}$ in $\left(
W^{*}(G),E\right) $ and hence $L_{w_{1}}$ and $L_{w_{2}}$ are not free over $%
D_{G}.$ This contradict our assumption.

\strut

(Case II) Suppose that the finite paths $w_{1},w_{2}$ are non-loop finite
paths in $loop^{c}(G)$ and assume that they are not diagram-distinct. Since
they are not diagram-distinct, they are identically equal. Therefore, they
are not free over $D_{G}$ in $\left( W^{*}(G),E\right) .$

\strut

(Case III) Let $w_{1}\in loop(G)$ and $w_{2}\in loop^{c}(G).$ They are
always diagram-distinct.

\strut

Let $L_{w_{1}}$ and $L_{w_{2}}$ are free over $D_{G}$ and assume that $w_{1}$
and $w_{2}$ are not diagram-distinct. Then $L_{w_{1}}$ and $L_{w_{2}}$ are
not free over $D_{G}$, by the Case I, II and III. So, this contradict our
assumption.
\end{proof}

\strut

The previous theorem characterize the $D_{G}$-freeness of two partial
isometries $L_{w_{1}}$ and $L_{w_{2}},$ where $w_{1},w_{2}\in FP(G).$ This
characterization shows us that the diagram-distinctness of finite paths
determine the $D_{G}$-freeness of corresponding creation operators.

\strut

Let $a$ and $b$ be the given $D_{G}$-valued random variables in (3.0). We
can get the necessary condition for the $D_{G}$-freeness of $a$ and $b,$ in
terms of their supports. \strut \strut \strut Recall that we say that the
two subsets $X_{1}$ and $X_{2}$ of $FP(G)$ are said to be diagram-distinct
if $x_{1}$ and $x_{2}$ are diagram-distinct, for all pairs $(x_{1},x_{2})$ $%
\in $ $X_{1}$ $\times $ $X_{2}.$

\strut \strut

\begin{theorem}
Let $a,b\in \left( W^{*}(G),E\right) $ be $D_{G}$-valued random variables
with their supports $\mathbb{F}^{+}(G:a)$ and $\mathbb{F}^{+}(G:b).$ The $%
D_{G}$-valued random variables $a$ and $b$ are free over $D_{G}$ in $\left(
W^{*}(G),E\right) $ if $FP(G:a_{1})$ and $FP(G:a_{2})$ are diagram-distinct.
\end{theorem}

\strut

\begin{proof}
For convenience, let's denote $a$ and $b$ by $a_{1}$ and $a_{2},$
respectively. Assume that the supports of $a_{1}$ and $a_{2},$ $\mathbb{F}%
^{+}(G:a_{1})$ and $\mathbb{F}^{+}(G:a_{2})$ are diagram-distinct. Then by
the previous $D_{G}$-freeness characterization,

\strut

$\ \ \ \underset{l\in FP(G:a_{1}),\,u_{l}\in \{1,*\}}{\sum }%
p_{l}^{(1)}L_{l}^{u_{l}}$ \ and \ $\underset{l\in FP(G:a_{2}),\,u_{l}\in
\{1,*\}}{\sum }p_{l}^{(2)}L_{l}^{u_{l}}$

\strut

are free over $D_{G}$ in $\left( W^{*}(G),E\right) .$ Indeed, since $%
FP(G:a_{1})$ and $FP(G:a_{2})$ are diagram-distinct, all summands $L_{w_{1}}$%
's of $a_{1}$ and $L_{w_{2}}$'s of $a_{2}$ are free over $D_{G}$ in $\left(
W^{*}(G),E\right) .$ Therefore, $a_{1}$ and $a_{2}$ are free over $D_{G}$ in 
$\left( W^{*}(G),E\right) .$\strut \strut \strut
\end{proof}

\strut \strut

\strut

\strut \strut

\section{$D_{G}$-valued Semicircular Elements}

\bigskip

\bigskip

Throughout this chapter, we will consider the $D_{G}$-valued
semicircularity. Let $B$ be a von Neumann algebra and $A$, a $W^{\ast }$%
-algebra over $B$ and let $E:A\rightarrow B$ be a conditional expectation.
Then $(A,E)$ be an amalgamated $W^{\ast }$-probability space over $B.$ We
say that the $D_{G}$-valued random variable $x\in \left( A,E\right) $ is a $%
B $-valued semicircular element if it is self-adjoint and the only
nonvanishing $B$-valued cumulants is the second one. i.e., $x\in (A,E)$ is $%
B $-semicircular if $x$ is self-adjoint in $A$ and

\bigskip

\begin{center}
$k_{n}\left( x,...,x\right) =\left\{ 
\begin{array}{ccc}
k_{2}(x,x)\neq 0_{B} &  & \text{if }n=2 \\ 
&  &  \\ 
0_{B} &  & \text{otherwise.}
\end{array}
\right. $
\end{center}

\bigskip

Let $G$ be a directed graph with $V(G)=\{v\}$ and $E(G)=\{e=vev\}.$ i.e., $G$
is a one-vertex graph with one loop edge. Canonically, we can construct the
graph $W^{\ast }$-algebra $W^{\ast }(G)$ and its diagonal subalgebra $%
D_{G}\simeq \Delta _{1}=\mathbb{C}.$ So, the canonical conditional
expectation $E$ is a linear functional. Moreover, it is a trace on $W^{\ast
}(G).$ Notice that, by Voiculescu,

\bigskip

\begin{center}
$W^{\ast }(G)\simeq L(\mathbb{Z}),$
\end{center}

\strut \strut

where $L(\mathbb{Z})$ is the free group factor. Moreover, the random
variable $L_{e}+L_{e}^{*}$ is the Voiculescu's semicircular element.
Therefore, the only nonvanishing cumulants of $L_{e}+L_{e}^{*}$ is the
second one. First, we will consider the following combinatorial fact. This
is crucial to consider the $D_{G}$-semicircularity on $\left(
W^{*}(G),E\right) .$\bigskip

\strut

\begin{lemma}
Let $G$ and $e$ be given as above. Then

\strut 

$\ \ \ \ \ \ \ \ \ \ \ \ \ \ \ \underset{L\in LP_{2}^{*}}{\sum }\mu
_{e,e}^{L(u_{1},u_{2})}=2$

and

$\ \ \ \ \ \ \ \underset{L\in LP_{n}^{*}}{\sum }\mu
_{e,...,e}^{L(u_{1},...,u_{n})}=0,$ $\forall n\in 2\mathbb{N}\setminus \{2\}.
$
\end{lemma}

\bigskip

\begin{proof}
Define $a=L_{e}+L_{e}^{\ast }.$ Then it is a semicircular element, in the
sense of Voiculescu. So, the only nonvanishing $D_{G}$-valued cumulant of $a$
is the second one. By the previous lemma, we have that

\strut

$\ \ \ \ \ \ \ \ \ \ \ \ \ k_{n}\left( a,...,a\right) =\underset{L\in
LP_{n}^{\ast }}{\sum }\mu _{e,...,e}^{L(u_{1},...,u_{n})}\cdot L_{v}.$

\strut

Suppose that $n=2.$ Then

\strut

$\ \ \ \ \ \ \ k_{2}(a,a)=\mu _{e,e}^{1,\ast }\cdot L_{v}+\mu _{e,e}^{\ast
,1}\cdot L_{v}=\mu _{e,e}^{1,\ast }+\mu _{e,e}^{\ast ,1},$

\strut

since $L_{v}=1_{D_{G}}=1_{\mathbb{C}}=1.$ Notice that $l_{e}l_{e^{-1}}$ and $%
l_{e^{-1}}l_{e}$ have their lattice paths

\strut

$\ \ \ \ \ \ \ \ \ \ \ \ \ \ \ \ _{\ast }\nearrow \searrow $ \ \ \ \ \ \ \ \
and \ \ \ \ \ \ \ $^{\ast }\searrow \nearrow ,$

\strut

respectively. Therefore, $C_{e,e}^{1,\ast }=C_{e,e}^{\ast ,1}=\{1_{2}\}.$
Therefore,

\strut

$\ \ \ \ \ \ \ \ \ \ \ \ \ \ \mu _{e,e}^{1,\ast }=\mu _{e,e}^{\ast ,1}=\mu
(1_{2},1_{2})=1.$

\strut

Therefore, $k_{2}(a,a)=2.$ Equivalently,

\strut

$\ \ \ \ \ \ \ \ \ \ \ \ \ \ \ \ \ \ \ \ \ \ \ \ \ \ \underset{L\in
LP_{2}^{\ast }}{\sum }\mu _{e,e}^{L(u_{1},u_{2})}=2,$

\strut

since $L_{v}=1_{D_{G}}=1_{\mathbb{C}}=1.$ Now, let $2<n\in 2\mathbb{N}.$ Then

\strut

$\ \ \ \ \ \ \ k_{n}\left( a,...,a\right) =\underset{L\in LP_{n}^{\ast }}{%
\sum }\mu _{e,...,e}^{L(u_{1},...,u_{n})}\cdot L_{v}=\underset{L\in
LP_{n}^{\ast }}{\sum }\mu _{e,...,e}^{L(u_{1},...,u_{n})}$

\strut

since $L_{v}=1$

\strut

$\ \ \ \ \ \ \ \ \ \ \ \ \ \ \ \ \ \ \ \ \ \ \ \ \ \ =0,$

\strut

by the semicircularity of $L_{e}+L_{e}^{\ast }.$ Therefore,

\strut

$\ \ \ \ \ \ \ \ \ \ \ \ \ \ \ \ \underset{L\in LP_{n}^{*}}{\sum }\mu
_{e,...,e}^{L(u_{1},...,u_{n})}=0,$ $\forall n\in 2\mathbb{N}~\setminus
~\{2\}.\strut $
\end{proof}

\bigskip \strut

By the previous lemma, we can determine the $D_{G}$\strut -semicircular
elements in $\left( W^{*}(G),E\right) .$

\strut

\begin{theorem}
Let $G$ be a countable directed graph and let $(W^{*}(G),E)$ be the graph $%
W^{*}$-probability space over the diagonal subalgebra $D_{G}.$ Let $w=vwv\in
loop(G),$ with $v\in V(G).$ Then $L_{w}+L_{w}^{*}$ is a $D_{G}$-valued
semicircular elements, with

\strut \strut 

$\ \ \ \ \ \ \ \ \ \ \ \ k_{2}(L_{w}+L_{w}^{*},L_{w}+L_{w}^{*})=2L_{v}.$
\end{theorem}

\bigskip

\begin{proof}
Let $w=vwv\in loop(G),$ with $v\in V(G).$ Define a $D_{G}$-valued random
variable $a=L_{w}+L_{w}^{\ast }.$ Then, clearly, it is self-adjoint in $%
W^{\ast }(G).$ It suffices to show that it has only second nonvanishing $%
D_{G}$-valued cumulants.

\strut

$\ \ \ k_{n}\left( a,...,a\right) =k_{n}\left( L_{w}+L_{w}^{\ast
},...,L_{w}+L_{w}^{\ast }\right) $

\strut

$\ \ \ \ \ \ \ =\underset{(u_{1},...,u_{n})\in \{1,\ast \}^{n}}{\sum }%
k_{n}\left( L_{w}^{u_{1}},...,L_{w}^{u_{n}}\right) $

\strut

by the bimodule map property of $k_{n}$

\strut

$\ \ \ \ \ \ \ =\underset{l_{w^{t_{1}}}...l_{w^{t_{n}}}:\ast \text{%
-axis-property}}{\sum }k_{n}\left( L_{w}^{u_{1}},...,L_{w}^{u_{n}}\right) $

\strut

where $l_{w^{t_{1}}}...l_{w^{t_{n}}}\in \left( A_{G^{\symbol{94}%
}}~/~R,~E\right) $ is the corresponding element of $%
L_{w}^{u_{1}}...L_{w}^{u_{n}}$

\strut

$\ \ \ \ \ \ \ =\underset{L\in LP_{n}^{\ast }}{\sum }\mu
_{w,...,w}^{L(u_{1},...,u_{n})}\cdot L_{v}.$

\strut

But by the previous theorem,

\strut

$\ \ \ \ \ \ \ \ \ \ \ \underset{L\in LP_{n}^{\ast }}{\sum }\mu
_{w,...,w}^{L(u_{1},...,u_{n})}=\left\{ 
\begin{array}{ccc}
2 &  & \text{if }n=2 \\ 
&  &  \\ 
0 &  & \text{otherwise.}
\end{array}
\right. $

\strut

Therefore,

\strut

$\ \ \ \ \ \ \ \ \ \ \ \ \ k_{n}(a,...,a)=\left\{ 
\begin{array}{ccc}
2L_{v} &  & \text{if }n=2 \\ 
&  &  \\ 
0_{D_{G}} &  & \text{otherwise.}
\end{array}
\right. $

\strut

and hence $a=L_{w}+L_{w}^{\ast },$ $w\in loop(G),$ is a $D_{G}$-semicircular
element in $(W^{\ast }(G),E).$ Notice that if $w$ is a loop (as a finite
path), then $%
C_{w,...,w}^{L(u_{1},...,u_{n})}=C_{e,...,e}^{L(u_{1},...,u_{n})},$ where $e$
is the given at the beginning of this chapter. Therefore, we have the
previous formuli.
\end{proof}

\bigskip \strut

So, we can conclude that all $D_{G}$-valued random variables having the
forms of $pL_{w}+pL_{w}^{\ast }$ \ ($p\in \mathbb{R},$ $w\in loop(G)$) are $%
D_{G}$-semicircular elements in $\left( W^{\ast }(G),E\right) .$

\bigskip

\begin{corollary}
Let $w_{j}=v_{j}w_{j}v_{j}\in loop(G),$ $v_{j}\in V(G),$ \ $j=1,...,N,$ such
that they are \textbf{mutually diagram-distinct} (Note that it is possible
that $v_{i}=v_{j},$ for any $i\neq j$ in $\{1,...,N\}$). Then the $D_{G}$%
-valued random variable

\strut 

$\ \ \ \ \ \ \ a=\sum_{j=1}^{N}\left(
p_{j}L_{w_{j}}+p_{j}L_{w_{j}}^{*}\right) \in \left( W^{*}(G),E\right) $

\strut 

is $D_{G}$-semicircular and

\strut 

$\ \ \ \ \ \ \ \ \ \ \ k_{2}\left( a,a\right) =\sum_{j=1}^{N}2p_{j}^{4}\cdot
L_{v_{j}}.$

$\square $
\end{corollary}

\bigskip \bigskip

\strut

\section{$D_{G}$-Even Elements}

\strut \bigskip \strut

\strut

Let $G$ be a directed graph with $V(G)=\{v_{1},v_{2}\}$ and $%
E(G)=\{e=v_{1}ev_{2}\}.$ We can construct the graph $W^{\ast }$-algebra $%
W^{\ast }(G)$ and its diagonal subalgebra $D_{G}=\Delta _{2}.$ Trivially $%
L_{v_{1}}+L_{v_{2}}=1_{D_{G}}.$ Define a $D_{G}$-valued random variable $%
a=L_{e}+L_{e}^{\ast }.$ Note that

$\bigskip $

\begin{center}
$FP(G)=E(G),$
\end{center}

\bigskip

since $e^{k}\notin \mathbb{F}^{+}(G),$ for all $k\in \mathbb{N}~\setminus
~\{1\}$. We cannot construct the more finite paths other than $e,$ itself.
Also, the support of this operator $a$ is

\bigskip

\begin{center}
$\mathbb{F}^{+}(G:a)=FP(G:a)=FP_{\ast }(G:a)=E(G).$
\end{center}

\bigskip

\bigskip Similar to the previous chapter, we will observe the following
combinatorial fact ;

\strut

\begin{lemma}
Let $G$ and $e$ be given as above. Then

\strut 

\ \ \ $\ \ \ \ \ \ \ \underset{L\in LP_{n}^{*}}{\sum }\mu
_{e,...,e}^{L(u_{1},...,u_{n})}=2\mu _{n}\in \mathbb{R},$

\strut 

where $\mu _{n}=\mu _{e,e,...,e,e}^{1,*,...,1,*}=\mu
_{e,e,...,e,e}^{*,1,...,*,1},$ for all $n\in 2\mathbb{N}.$
\end{lemma}

\bigskip

\begin{proof}
Since $e=v_{1}ev_{2},$ with $v_{1}\neq v_{2}\in V(G),$ $e^{k}\notin \mathbb{F%
}^{+}(G)$ (i.e., it is not admissible), whenever $n>1.$ Therefore, for any
even $n,$

\strut

$\ \ \ \ \ \ \ \ \ \underset{L\in LP_{n}^{\ast }}{\sum }\mu
_{e,...,e}^{L(u_{1},...,u_{n})}=\mu _{e,e,...,e,e}^{1,\ast ,...,1,\ast }+\mu
_{e,e,...,e,e}^{\ast ,1,...,\ast ,1}.$

\strut

Indeed, if we let $(u_{1},...,u_{n})\in \{1,\ast \}^{n}$ and if it is not
alternating, then there exists at least one $j\in \{1,...,n\}$ such that $%
u_{j}=u_{j+1}$ in $\{1,\ast \}.$ This means that there should be a
consecutive increasing or decreasing words of $e$'s. But $e^{k}$ does not
exist in our graph $G,$ whenever $k>1.$ Therefore, for such $%
L(u_{1},...,u_{n})\in LP_{n}^{\ast },$

\strut

$\ \ \ \ \ \ \ \ \ \ \ \mu _{e,...,e}^{L(u_{1},...,u_{n})}=\underset{\theta
\in C_{e,...,e}^{L(u_{1},...,u_{n})}}{\sum }\mu (\theta ,1_{n})=0\in \mathbb{%
R}.$

\strut

Observe that $C_{e,e,...,e,e}^{1,\ast ,...,1,\ast }$ and $%
C_{e,e,...,e,e}^{\ast ,1,...,\ast ,1}$ have the same elements, by the
symmetry of the lattice paths of $l_{e}l_{e^{-1}}...l_{e}l_{e^{-1}}$ and $%
l_{e^{-1}}l_{e}...l_{e^{-1}}l_{e}.$ So,

\strut

$\ \ \ \ \ \ \ \ \ \ \ \ \ \ \ \ \ \mu _{e,e,...,e,e}^{1,\ast ,...,1,\ast
}=\mu _{e,e,...,e,e}^{\ast ,1,...,\ast ,1}.$

\strut

Denote this real value by $\mu _{n},$ for each $n\in 2\mathbb{N}.$ Then we
can get the above formula.
\end{proof}

\bigskip \strut

So, we have that ;

\bigskip

\begin{lemma}
Let $G$ and $e$ be given as above. Then the $D_{G}$-valued cumulants of $%
a=L_{e}+L_{e}^{*}$ is determined by ;

\strut 

(1) $k_{n}\left( a,...,a\right) =0_{D_{G}},$ whenever $n$ is odd.

\strut 

(2) $k_{n}(a,...,a)=\mu _{n}\cdot 1_{D_{G}},$ for all $n\in 2\mathbb{N},$
where

\strut 

$\ \ \ \ \ \ \ \ \ \ \ \ \ \ \ \ \ \ \ \ \ \ \ \ \ \ \ \ \ \ \ \ \mu
_{n}=\mu _{e,e,...,e,e}^{1,*,...1,*}=\mu _{e,e,...,e,e}^{*,1,...,*,1}.$
\end{lemma}

\bigskip

\begin{proof}
Suppose $n$ is odd. Then there is no lattice path having the $\ast $%
-axis-property. Therefore, (1) holds true. Now, assume that $n\in 2\mathbb{N}%
.$ Let $L\in LP_{n}^{\ast }.$ Since $e$ is a non-loop edge, $e^{k}\notin
FP(G)=E(G).$ Therefore, when we consider $l_{e^{t_{1}}}...l_{e^{t_{n}}}\in
\left( A_{G^{\symbol{94}}}~/~R,~E\right) ,$ the corresponding element of $%
L_{e}^{u_{1}}...L_{e}^{u_{n}}$ (where $L_{e}$ and $L_{e}^{\ast }$ are
summands of $a$), the lattice paths of it is

\strut

$\ \ \ \ \ \ \ _{\ast }\nearrow \searrow \nearrow \searrow ...\nearrow
\searrow $ \ \ or \ \ $^{\ast }\searrow \nearrow \searrow \nearrow ^{\cdot
\cdot \cdot }\searrow \nearrow ,$

\strut

denoted by $[rfrf...rf]$ and $[frfr...fr],$ respectively (Here $r$ stands
for the rising step and $f$ stands for the falling step), because if there
are consecutive rising steps or consecutive falling steps, then it
represents $l_{e}...l_{e}=0_{D_{G}}$ or $l_{e^{-1}}...l_{e^{-1}}.$
Therefore, if $n\in 2\mathbb{N},$ then

\strut

$\ k_{n}(a,...,a)$

\strut

$\ \ \ =\underset{(u_{1},...,u_{n})\in \{1,\ast \}^{n}}{\sum }k_{n}\left(
L_{e}^{u_{1}},...,L_{e}^{u_{n}}\right) $

\strut

$\ \ \ =\underset{L\in LP_{n}^{*}}{\sum }\mu
_{e,....,e}^{L(u_{1},...,u_{n})}\cdot E(L_{e}^{u_{1}}...L_{e}^{u_{n}})$

\strut

$\ \ \ =\mu
_{e,e,....,e,e}^{1,*,...,1,*}E(L_{e}L_{e}^{*}...L_{e}L_{e}^{*})+\mu
_{e,e,....,e,e}^{*,1,...,*,1}E\left( L_{e}^{*}L_{e}...L_{e}^{*}L_{e}\right) $

\strut

$\ \ \ =\mu _{n}\cdot L_{v_{1}}+\mu _{n}L_{v_{2}}=\mu
_{n}(L_{v_{1}}+L_{v_{2}})\ $

\strut

where $\mu _{n}=\mu _{e,e,....,e,e}^{1,\ast ,...,1,\ast }=\mu
_{e,e,....,e,e}^{\ast ,1,...,\ast ,1}$

\strut

$\ \ \ =\mu _{n}\cdot 1_{D_{G}}.$

\strut

Indeed, the lattice paths $[rf...rf]$ and $[fr...fr]$ in $LP_{n}^{\ast }$
induce the same sets $C_{e,e,....,e,e}^{1,\ast ,...,1,\ast }$ and $%
C_{e,e,....,e,e}^{\ast ,1,...,\ast ,1}.$ Thus

\strut

$\ \ \ \ \mu _{e,e,....,e,e}^{1,\ast ,...,1,\ast }=\underset{\theta \in
C_{e,e,....,e,e}^{1,\ast ,...,1,\ast }}{\sum }\mu (\theta ,1_{n})=\underset{%
\pi \in C_{e,e,....,e,e}^{\ast ,1,..,.\ast ,1}}{\sum }\mu (\pi ,1_{n})=\mu
_{e,e,....,e,e}^{\ast ,1,...,\ast ,1}.$

\strut
\end{proof}

\bigskip

Now, we will observe the general case.\bigskip \ Define $\mu _{n}\overset{def%
}{=}\mu _{e,e,....,e,e}^{1,*,...,1,*}=\mu _{e,e,....,e,e}^{*,1,...,*,1},$
where $e$ is given as before. Now, we will observe the general case when $e$
is a general finite path in an arbitrary countable directed graph $G.$

\strut \strut

\begin{definition}
Let $B$ be a von Neumann algebra and let $A$ be the $W^{*}$-algebra over $B.$
Let $E_{B}:A\rightarrow B$ be a conditional expectation and $\left(
A,E_{B}\right) $, the amalgamated $W^{*}$-probability space over $B.$ We say
that the $B$-valued random variable $a\in \left( A,E_{B}\right) $ is $B$%
-valued even (in short, $B$-even) if it is self-adjoint and if it has all
vanishing odd $B$-valued moments.
\end{definition}

\strut

Let $B$ be a von Neumann algebra and $A,$ a $W^{*}$-algebra over $B$ and let 
$\left( A,E_{B}\right) $ be the amalgamated $W^{*}$-probability space over $%
B.$ Recall that $B$-valued random variable $a\in \left( A,E_{B}\right) _{sa}$
is $B$-even if and only if all odd $B$-valued cumulants of $a$ vanish. By
the previous lemma, we can easily see that $D_{G}$-valued random variables $%
L_{w}+L_{w}^{*},$ for $w\in FP(G),$ are $D_{G}$-even, because if $%
w=v_{1}wv_{2},$ then

\strut

\begin{center}
$k_{n}\left( L_{w}+L_{w}^{*},...,L_{w}+L_{w}^{*}\right) =\left\{ 
\begin{array}{ll}
0_{D_{G}} & \text{if }n\text{ is odd} \\ 
&  \\ 
\mu _{n}\cdot \left( L_{v_{1}}+L_{v_{2}}\right) , & \text{if }n\text{ is even%
}
\end{array}
\right. $
\end{center}

\strut

for all $n\in \mathbb{N},$ where $\mu _{n}=\mu
_{w,w,...,w,w}^{1,*,...1,*}=\mu _{w,w,...,w,w}^{*,1,...,*,1}.$ Remark that

$\strut $

$%
C_{w,w,...,w,w}^{1,*,...,1,*}=C_{w,w,...,w,w}^{*,1,...,*,1}=C_{e,e,...,e,e}^{*,1,...,*,1}=C_{e,e,...,e,e}^{1,*,...,1,*}, 
$

\strut

where $e$ is the edge given at the beginning of this chapter. Based on it,
we can get that

\strut

\begin{theorem}
Let $w\in FP(G).$ Then the $D_{G}$-valued random variable $L_{w}+L_{w}^{*}$
is $D_{G}$-even. Moreover, if $w=v_{1}wv_{2}$ with $v_{1},v_{2}\in V(G)$
(possibly $v_{1}=v_{2}$), then

\strut 

$\ k_{n}\left( L_{w}+L_{w}^{*},...,L_{w}+L_{w}^{*}\right) =\left\{ 
\begin{array}{ll}
0_{D_{G}} & \text{if }n\text{ is odd} \\ 
&  \\ 
\mu _{n}\cdot \left( L_{v_{1}}+L_{v_{2}}\right) , & \text{if }n\text{ is even%
}
\end{array}
\right. $

\strut 

for all $n\in \mathbb{N},$ where $\mu _{n}=\mu
_{w,w,...,w,w}^{1,*,...1,*}=\mu _{w,w,...,w,w}^{*,1,...,*,1}.$ $\square $
\end{theorem}

\strut

\strut

\strut

\section{\strut $D_{G}$-valued R-diagonal Elements}

\bigskip

\bigskip

In this chapter, we will consider the $D_{G}$-valued R-diagonality on the
graph $W^{\ast }$-probability space $\left( W^{\ast }(G),E\right) $ over the
diagonal subalgebra $D_{G}.$ Recall that

\bigskip

\begin{definition}
Let $B$ be a von Neumann algebra and let $A$ be a $W^{*}$-algebra over $B.$
Suppose that we have a conditional expectation $E_{B}:A\rightarrow B$ and
hence $(A,E_{B})$ is the amalgamated $W^{*}$-probability space over $B.$ We
say that the $D_{G}$-valued random variable $x\in (A,E_{B})$ is a $D_{G}$%
-valued R-diagonal element if the only nonvanishing mixed cumulants of $x$
and $x^{*}$ are

\strut 

$\ \ \ \ \ \ \ \ \ \ \ \ \ k_{n}(b_{1}x,b_{2}x^{*},...,b_{n-1}x,b_{n}x^{*})$

and

$\ \ \ \ \ \ \ \ \ \ \ \ \ k_{n}(b_{1}^{\prime
}x^{*},b_{2}x,...,b_{n-1}^{\prime }x^{*},b_{n}^{\prime }x),$

\strut 

for all $n\in 2\mathbb{N},$ where $b_{1},b_{1}^{\prime
},...,b_{n},b_{n}^{\prime }\in B$ are arbitrary. If $n$ is odd, then
automatically the mixed cumulants vanish.
\end{definition}

\bigskip

We can show that $L_{w}$ is $D_{G}$-valued R-diagonal, whenever $w$ is a
finite path in $G.$ By the results in the proceeding two chapters, we can
get the following theorem ;

\strut

\begin{theorem}
Let $G$ be a countable directed graph and $w\in FP(G).$ Then the $D_{G}$%
-valued random variable $L_{w}$ and $L_{w}^{*}\in \left( W^{*}(G),E\right) $
are $D_{G}$-valued R-diagonal.
\end{theorem}

\bigskip

\begin{proof}
It suffices to show that the only nonvanishing mixed cumulants of $L_{w}$
and $L_{w}^{*}$ are alternating ones. i.e., the nonvanishing mixed cumulants
are

\strut

(6.1)$\ \ k_{2n}\left( L_{w}^{*},L_{w},...,L_{w}^{*},L_{w}\right) $ \ and \ $%
k_{2n}\left( L_{w},L_{w}^{*},...,L_{w},L_{w}^{*}\right) .$

\strut

Suppose that $w$ is a loop. Then by Lemma 4.1, we have that the only
nonvanishing mixed cumulants are

\strut

$\ \ \ \ \ \ \ \ \ k_{2}\left( L_{w}^{*},L_{w}\right) $ \ \ \ and \ \ \ $%
k_{2}\left( L_{w},L_{w}^{*}\right) .$

\strut

So, $L_{w}$ and $L_{w}^{*}$ are $D_{G}$-valued R-diagonal.

\strut

Now, assume that $w=v_{1}wv_{2}$ is a non-loop finite path, with $v_{1}\neq
v_{2}\in V(G).$ Then the nonvanishing mixed cumulants of $L_{w}$ and $%
L_{w}^{*}$ have the forms of (6.1). By Section 2.2, we can easily get that

\strut

$\ \ \ \ \ \ \ k_{2n}\left( L_{w}^{*},L_{w}...,L_{w}^{*},L_{w}\right) =\mu
_{w,w,...,w,w}^{*,1,...,*,1}\cdot L_{v_{1}}$

and

$\ \ \ \ \ \ \ k_{2n}\left( L_{w},L_{w}^{*},...,L_{w},L_{w}^{*}\right) =\mu
_{w,w,...,w,w}^{1,*,...,1,*}\cdot L_{v_{2}}.$

\strut \strut

Assume that there exists a nonvanishing mixed $D_{G}$-cumulant of $L_{w}$
and $L_{w}^{*}.$ i.e., assume that there exist $n\in \Bbb{N}$ and a $n$%
-tuple $(u_{1},...,u_{n})$ of $\{1,*\}$ such that $k_{n}\left(
L_{w}^{u_{1}},...,L_{w}^{u_{n}}\right) \neq 0_{D_{G}}.$ By Section 2.2, we
have that

\strut \strut

$\ \ \ \ \ \ \ k_{n}\left( L_{w}^{u_{1}},...,L_{w}^{u_{n}}\right) =\mu
_{w,...,w}^{u_{1},...,u_{n}}E(L_{w}^{u_{1}}...L_{w}^{u_{n}}).$

\strut

Notice that since $w$ is a non-loop finite path, there is no admissible
finite path $w^{k},$ for $k\in \Bbb{N\,}\setminus \,\{1\}$. So, if $%
(u_{1},...,u_{n})$ is not alternating, then there exists at least one $j$ in 
$\{1,...,n-1\}$ such that $u_{j}=u_{j+1}.$ Since $%
L_{w}^{u_{j}}L_{w}^{u_{j+1}}=L_{w}^{2}$ or $L_{w}^{*\,\,2},$ the $D_{G}$%
-valued random variable $L_{w}^{u_{1}}...L_{w}^{u_{n}}$ does not have the $*$%
-axis-property and hence

\strut

$\ \ \ \ \ \ \
E(L_{w}^{u_{1}}...L_{w}^{u_{n}})=0_{D_{G}}=k_{n}(L_{w}^{u_{1}},...,L_{w}^{u_{n}}). 
$

\strut

This contradict our assumption. So, $L_{w}$ and $L_{w}^{*}$ are $D_{G}$%
-valued R-diagonal.\strut
\end{proof}

\bigskip \strut \strut

The above theorem shows us that all generators of $W^{*}(G)$ generated by
finite paths in $FP(G)$ are $D_{G}$-valued R-diagonal.

\strut \strut

\strut \strut

\section{Generating Operators}

\strut \strut

\strut

\strut

In this chapter, as examples, we will compute the trivial $D_{G}$-valued
moments and cumulants of the generating operator $T$ of the graph $W^{*}$%
-algebra $W^{*}(G).$ Let $G$ be a countable directed graph and let $\left(
W^{*}(G),E\right) $ be the graph $W^{*}$-probability space over its diagonal
subalgebra. Let $a\in \left( W^{*}(G),E\right) $ be a $D_{G}$-valued random
variable. Recall that the trivial $D_{G}$-valued $n$-th moments and
cumulants of $a$ are defined by

\strut

\begin{center}
$E(a^{n})$ \ \ and \ $k_{n}\left( \underset{n\text{-times}}{\underbrace{%
a,.......,a}}\right) .$
\end{center}

\strut

In this chapter, we will deal with the following special $D_{G}$-valued
random variable ;

\strut

\begin{definition}
Define an operator $T$ in $W^{*}(G)$ by

\strut 

(5.1)$\strut $ $\ \ \ \ \ \ \ \ \ \ \ \ \ T=\underset{e\in E(G)}{\sum }%
\left( L_{e}+L_{e}^{*}\right) .$

\strut 

We will call $T$ the generating operator of $W^{*}(G).$ The self-adjoint
operators $L_{e}+L_{e}^{*},$ for $e\in E(G),$ are called the block operators
of $T.$
\end{definition}

\strut

\begin{example}
Let $G$ be a one-vertex directed graph with $N$-edges. i.e.,

\strut 

$\ \ \ \ V(G)=\{v\}$ \ and \ $E(G)=\{e_{j}=ve_{j}v:j=1,...,N\}.$

\strut 

Then the graph $W^{*}$\strut -algebra $W^{*}(G)$ satisfies that

\strut 

(5.2)$\ \ \ \ W^{*}(G)=D_{G}*_{D_{G}}\left( \underset{j=1}{\overset{N}{%
\,*_{D_{G}}}}\left( W^{*}(\{L_{e_{j}}\},D_{G})\right) \right) ,$

\strut 

by Chapter 4. Notice that $D_{G}=\Delta _{1}=\mathbb{C}.$ Therefore, the
formula (5.2) is rewritten by

\strut 

(5.4) \ \ \ $\ \ \ \ \ \ \ W^{*}(G)=\,\underset{j=1}{\overset{N}{*}}\left(
W^{*}(\{L_{e_{j}}\})\right) ,$

\strut 

where $*$ means the usual (scalar-valued) free product of Voiculescu. Also
notice that $1_{D_{G}}=L_{v}=1\in \mathbb{C}$ and

\strut 

(5.5) $\ \ \ L_{e_{j}}^{*}L_{e_{j}}=L_{v}=1=L_{e_{j}}L_{e_{j}}^{*},$ for all 
$j=1,...,N.$

\strut 

This shows that $L_{e_{j}}$'s are unitary in $W^{*}(G),$ for all $j=1,...,N.$
Now, define the generating operator $T=\sum_{j=1}^{N}\left(
L_{e_{j}}+L_{e_{j}}^{*}\right) $ of $W^{*}(G).$ It is easy to see that each
block operator $L_{e_{j}}+L_{e_{j}}^{*}$ is semicircular, by Voiculescu, for
all $j=1,...,N.$ (Remember the construction of creation operators $L_{e_{j}}$%
's and see [9].) Futhermore, by Chapter 3, we can get that all blocks $%
L_{e_{j}}+L_{e_{j}}^{*}$'s are free from each other in the graph $W^{*}$%
-probability space $\left( W^{*}(G),E\right) .$

\strut 

By (5.4), the canonical conditional expectation $E:W^{*}(G)\rightarrow D_{G}$
is the faithful linear functional. Moreover, by (5.5), this linear
functional $E$ is a trace in the sense that $E(ab)=E(ba),$ for all $a,b\in
W^{*}(G).$ From now, to emphasize that $E$ is a trace, we will denote $E$ by 
$tr.$

\strut 

Let's compute the $n$-th cumulant of $T$ ;

\strut 

\strut (5.6)

$\ 
\begin{array}{ll}
k_{n}\left( T,...,T\right)  & =k_{n}\left( \sum_{j=1}^{N}\left(
L_{e_{j}}+L_{e_{j}}^{*}\right) ,...,\sum_{j=1}^{N}\left(
L_{e_{j}}+L_{e_{j}}^{*}\right) \right)  \\ 
&  \\ 
& =\sum_{j=1}^{N}k_{n}\left(
(L_{e_{j}}+L_{e_{j}}^{*}),...,(L_{e_{j}}+L_{e_{j}}^{*})\right) ,
\end{array}
$ $\ $

\strut \strut 

by the mutual freeness of $\{L_{e_{j}},L_{e_{j}}^{*}\}$'s on $\left(
W^{*}(G),tr\right) $, for $j=1,...,N.$ Observe that

\strut 

$\ \ \ k_{n}\left(
(L_{e_{j}}+L_{e_{j}}^{*}),...,(L_{e_{j}}+L_{e_{j}}^{*})\right) $

\strut 

(5.7) $\ \ \ =\left\{ 
\begin{array}{ll}
k_{2}\left( (L_{e_{j}}+L_{e_{j}}^{*}),(L_{e_{j}}+L_{e_{j}}^{*})\right)  & 
\text{if }n=2 \\ 
&  \\ 
0 & \text{otherwise,}
\end{array}
\right. $

\strut 

by the semicircularity of $L_{e_{j}}+L_{e_{j}}^{*},$ for $j=1,...,N.$ By
(5.7), the formula (5.6) is

\strut 

$\ \ \ k_{n}\left( T,...,T\right) $

\strut 

(5.8) $\ \ \ \ =\left\{ 
\begin{array}{ll}
\sum_{j=1}^{N}k_{2}\left(
L_{e_{j}}+L_{e_{j}}^{*},L_{e_{j}}+L_{e_{j}}^{*}\right)  & \text{if }n=2 \\ 
&  \\ 
0 & \text{otherwise}
\end{array}
\right. $

$\strut \ \ \ \ $

\strut Now, observe $k_{2}\left(
L_{e_{j}}+L_{e_{j}}^{*},L_{e_{j}}+L_{e_{j}}^{*}\right) $ ;

\strut 

$\ k_{2}\left( L_{e_{j}}+L_{e_{j}}^{*},L_{e_{j}}+L_{e_{j}}^{*}\right) $

\strut 

$\ \ \ =k_{2}\left( L_{e_{j}},L_{e_{j}}\right) +k_{2}\left(
L_{e_{j}},L_{e_{j}}^{*}\right) +k_{2}\left( L_{e_{j}}^{*},L_{e_{j}}\right)
+k_{2}\left( L_{e_{j}}^{*},L_{e_{j}}^{*}\right) $

\strut 

$\ \ \ =0+k_{2}\left( L_{e_{j}},L_{e_{j}}^{*}\right) +k_{2}\left(
L_{e_{j}}^{*},L_{e_{j}}\right) +0$

\strut 

by Section 2.1

\strut 

$\ \ \ =tr\left( L_{e_{j}}L_{e_{j}}^{*}\right) +tr\left(
L_{e_{j}}^{*}L_{e_{j}}\right) =2\cdot tr\left( L_{e_{j}}^{*}L_{e_{j}}\right) 
$

\strut 

since $tr$ is a trace

\strut 

$\ \ \ =2\cdot L_{v}=2,$

\strut 

for $j=1,...,N,$ by Section 2.1 and 2.2. So, we can get that

\strut 

(5.9) \ $k_{n}\left( T,...,T\right) =\left\{ 
\begin{array}{lll}
2N &  & \text{if }n=2 \\ 
&  &  \\ 
0 &  & \text{otherwise.}
\end{array}
\right. $

\strut 

Now, we can compute the trivial moments of $T,$ via the M\"{o}bius inversion.

\strut 

$\ tr\left( T^{n}\right) =\underset{\pi \in NC(n)}{\sum }k_{\pi }\left(
a,...,a\right) $

\strut 

where $k_{\pi }(a,...,a)=\underset{V\in \pi }{\prod }k_{\left| V\right|
}\left( \underset{\left| V\right| \text{-times}}{\underbrace{a,.......,a}}%
\right) ,$ for each $\pi \in NC(n),$ by Nica and Speicher (See [1] and [17])

\strut 

$\ \ \ \ \ \ =\underset{\pi \in NC_{2}(n)}{\sum }k_{\pi }(a,...,a)=\underset{%
\pi \in NC_{2}(n)}{\sum }\,\left( \underset{V\in \pi }{\prod }k_{\left|
V\right| }\left( a,a\right) \right) $

\strut 

where $NC_{2}(n)=\{\pi \in NC(n):V\in \pi \Leftrightarrow \left| V\right|
=2\}$ is the collection of all noncrossing pairings

\strut 

(5.10) $\ $

$\ \ \ \ \ \ =\underset{\pi \in NC_{2}(n)}{\sum }\,\left( \underset{V\in \pi 
}{\prod }2N\right) =\underset{\pi \in NC_{2}(n)}{\sum }\left( 2N\right)
^{\left| \pi \right| },$

\strut 

where $\left| \pi \right| \overset{def}{=}$ the number of blocks in $\pi .$
Notice that the above formula (5.10) shows us that the $n$ should be even,
because $NC_{2}(n)$ is nonempty when $n$ is even. Therefor,

\strut 

(5.11) \ \ \ $tr\left( T^{n}\right) =\left\{ 
\begin{array}{ll}
\underset{\pi \in NC_{2}(n)}{\sum }\left( 2N\right) ^{\left| \pi \right| } & 
\text{if }n\text{ is even} \\ 
&  \\ 
0 & \text{if }n\text{ is odd.}
\end{array}
\right. $

\strut 

Also, notice that if $\pi \in NC_{2}(n),$ then $\left| \pi \right| =\frac{n}{%
2},$ for all even number $n\in \mathbb{N}.$ So,

\strut 

\ \ \ $tr\left( T^{n}\right) =\left\{ 
\begin{array}{ll}
\left| NC_{2}(n)\right| \cdot \left( 2N\right) ^{\frac{n}{2}} & \text{if }n%
\text{ is even} \\ 
&  \\ 
0 & \text{if }n\text{ is odd}
\end{array}
\right. $

\strut 

(5.12) $\ \ \ =\left\{ 
\begin{array}{ll}
\left( 2N\right) ^{\frac{n}{2}}\cdot c_{\frac{n}{2}} & \text{if }n\text{ is
even} \\ 
&  \\ 
0 & \text{if }n\text{ is odd,}
\end{array}
\right. $

\strut 

where $c_{k}=\frac{1}{k+1}\left( 
\begin{array}{l}
2k \\ 
\,\,k
\end{array}
\right) $ is the $k$-th Catalan number, for all $k\in \mathbb{N}.$ Remember
that

\strut 

$\ \ \ \ \ \ \ \ \ \left| NC(k)\right| =\left| NC_{2}(2k)\right| =c_{k},$
for all $k\in \mathbb{N}.$

\strut 

Therefore, by (5.9) and (5.12), we can compute the moments and cumulants of
the generating operator $T$ of $\left( W^{*}(G),tr\right) $ ;

\strut 

$\ \ \ \ \ \ \ \ \ tr\left( T^{n}\right) =\left\{ 
\begin{array}{ll}
\left( 2N\right) ^{\frac{n}{2}}\cdot c_{\frac{n}{2}} & \text{if }n\text{ is
even} \\ 
&  \\ 
0 & \text{if }n\text{ is odd,}
\end{array}
\right. $

and

$\ \ \ \ \ \ \ \ \ k_{n}\left( T,...,T\right) =\left\{ 
\begin{array}{lll}
2N &  & \text{if }n=2 \\ 
&  &  \\ 
0 &  & \text{otherwise.}
\end{array}
\right. $

\strut 
\end{example}

\strut

\begin{example}
Let $N\in \mathbb{N}$ and let $G$ be the circulant graph with

\strut 

$\ \ \ \ \ \ \ \ \ \ \ \ \ V(G)=\{v_{1},...,v_{N}\}$

and

$\ \ \ \ \ \ \ \ \ \ \ \ \ E(G)=\{e_{1},...,e_{N}\}$

with

\strut 

$\ \ e_{j}=v_{j}e_{j}v_{j+1}$, for $j=1,...,N-1,$ and $e_{N}=v_{N}e_{N}v_{1}.
$

\strut 

Define the generating operator $T=\sum_{j=1}^{N}\left(
L_{e_{j}}+L_{e_{j}}^{*}\right) $ of the graph $W^{*}$-algebra $W^{*}(G).$ In
this case, we can get the diagonal subalgebra $D_{G}$ of $W^{*}(G),$ as a
von Neumann algebra which is isomorphic to $\Delta _{N},$ where $\Delta _{N}$
is a subalgebra of the matricial algebra $M_{N}(\mathbb{C}).$ Define the
canonical conditional expectation $E:W^{*}(G)\rightarrow D_{G}.$ Then we can
compute the trivial $n$-th $D_{G}$-valued cumulant of the operator $T,$ by
regarding it as a $D_{G}$-valued random variable in the graph $W^{*}$%
-probability space $\left( W^{*}(G),E\right) $ over $D_{G}=\Delta _{N}.$
Notice that each block $L_{e_{j}}+L_{e_{j}}^{*}$'s are free from each other
over $D_{G}$ in $\left( W^{*}(G),E\right) ,$ by the diagram-distinctness of $%
e_{j}$'s, for $j=1,...,N.$

\strut 

Fix $n\in \mathbb{N}.$ Then

\strut 

$\ k_{n}\left( \underset{n\text{-times}}{\underbrace{T,.......,T}}\right)
=k_{n}\left( \sum_{j=1}^{N}\left( L_{e_{j}}+L_{e_{j}}^{*}\right)
,...,\sum_{j=1}^{N}\left( L_{e_{j}}+L_{e_{j}}^{*}\right) \right) $

\strut 

$\ \ \ \ \ \ \ \ \ \ \ =\sum_{j=1}^{N}k_{n}\left(
(L_{e_{j}}+L_{e_{j}}^{*}),...,(L_{e_{j}}+L_{e_{j}}^{*})\right) $

\strut 

by the mutual $D_{G}$-freeness of $\{L_{e_{j}},L_{e_{j}}^{*}\}$'s, for $%
j=1,...,N$

\strut 

(5.13) $\ =\sum_{j=1}^{N}\underset{(u_{1},...,u_{n})\in \{1,*\}}{\sum }%
k_{n}\left( L_{e_{j}}^{u_{1}},...,L_{e_{j}}^{u_{n}}\right) .$

\strut 

Recall that, by Section 2.2, we can get that

\strut 

(5.14) $\ \ k_{n}\left( L_{e_{j}}^{u_{1}},...,L_{e_{j}}^{u_{n}}\right) =\mu
_{e_{j},...,e_{j}}^{u_{1},...,u_{n}}\cdot \Pr oj\left(
L_{e_{j}}^{u_{1}}...L_{e_{j}}^{u_{n}}\right) ,$

\strut 

where $\mu _{e_{j},...,e_{j}}^{u_{1},...,u_{n}}=\underset{\pi \in
C_{e_{j},...,e_{j}}^{u_{1},...,u_{n}}}{\sum }\mu (\pi ,1_{n}).$

\strut 

Observe that since $e_{j}$'s are non-loop edges, $e_{j}^{k}\notin \mathbb{F}%
^{+}(G),$ for all $k\in \mathbb{N}\,\setminus \,\{1\}$, for $j=1,...,N.$ In
other words, such $e_{j}^{k}$ is not admissible. So, if $(u_{1},...,u_{n})$
is not alternating, in the sense that $(u_{1},...,u_{n})=(1,*,...,1,*)$ or $%
(*,1,...,*,1),$ then $\Pr oj\left(
L_{e_{j}}^{u_{1}}...L_{e_{j}}^{u_{n}}\right) =0_{D_{G}}.$ For instance, $%
E\left( L_{e_{j}}^{*}L_{e_{j}}L_{e_{j}}^{*}\right) =0_{D_{G}}$ or $E\left(
L_{e_{j}}^{2}L_{e_{j}}^{*}L_{e_{j}}\right) =0_{D_{G}},$ by Section 2.1.
Therefore, the only nonvanishing case is either

\strut 

$\ \ \ k_{n}\left(
L_{e_{j}},L_{e_{j}}^{*},...,L_{e_{j}},L_{e_{j}}^{*}\right) $ \ or \ $%
k_{n}\left( L_{e_{j}}^{*},L_{e_{j}},...,L_{e_{j}}^{*},L_{e_{j}}\right) ,$

\strut 

where $n$ is even. Notice that

\strut 

(5.15) $\ \ \ \ \ \ \ \mu _{e_{j},e_{j},...,e_{j},e_{j}}^{1,*,...,1,*}=\mu
_{e_{j},e_{j},...,e_{j},e_{j}}^{*,1,...,*,1},$

\strut 

because $%
C_{e_{j},e_{j},...,e_{j},e_{j}}^{1,*,...,1,*}=C_{e_{j},e_{j},...,e_{j},e_{j}}^{*,1,...,*,1},
$ for all $j=1,...,N.$ Moreover, since $%
C_{e_{j},e_{j},...,e_{j},e_{j}}^{1,*,...,1,*}$ $=$ $%
C_{e_{k},e_{k},...,e_{k},e_{k}}^{1,*,...,1,*},$ for all $j\neq k$ in $%
\{1,...,N\},$

\strut 

(5.16) $\ \ \ \ \ \ \ \mu _{e_{j},e_{j},...,e_{j},e_{j}}^{1,*,...,1,*}=\mu
_{e_{k},e_{k},...,e_{k},e_{k}}^{1,*,...,1,*},$

\strut 

for all $j,k\in \{1,...,N\}.$ Let's denote $\mu
_{e_{j},e_{j},...,e_{j},e_{j}}^{1,*,...,1,*}$ by $\mu _{n},$ for all $%
j=1,...,N.$ Then, by (5.14), we have that

\strut 

\strut (5.17)

$k_{n}\left( L_{e_{j}}^{u_{1}},...,L_{e_{j}}^{u_{n}}\right) =\left\{ 
\begin{array}{ll}
\mu _{n}L_{v_{j}} & \text{if }(u_{1},...,u_{n})=(1,*,...,1,*) \\ 
&  \\ 
\mu _{n}L_{v_{j+1}} & \text{if }(u_{1},...,u_{n})=(*,1,...,*,1) \\ 
&  \\ 
0_{D_{G}} & \text{otherwise,}
\end{array}
\right. $

\strut 

for all $j=1,...,N,$ where $L_{v_{N+1}}$ means $L_{v_{1}}.$ So, by (5.13)
and (5.17), we can get that

\strut 

$\ k_{n}\left( T,...,T\right) $

\strut 

$\ =\sum_{j=1}^{N}\left( k_{n}\left(
L_{e_{j}},L_{e_{j}}^{*}...,L_{e_{j}},L_{e_{j}}^{*}\right) +k_{n}\left(
L_{e_{j}}^{*},L_{e_{j}},...,L_{e_{j}}^{*},L_{e_{j}}\right) \right) $

\strut \strut 

$\ =\sum_{j=1}^{N}\left( \mu _{n}L_{v_{j}}+\mu _{n}L_{v_{j+1}}\right)
=\sum_{j=1}^{N}\mu _{n}\left( L_{v_{j}}+L_{v_{j+1}}\right) $

\strut 

where $L_{v_{N+1}}$ means $L_{v_{1}}$, for all $n\in 2\mathbb{N}.$ Therefore,

\strut 

$\ \ k_{n}\left( T,...,T\right) =\left\{ 
\begin{array}{ll}
\sum_{j=1}^{N}\mu _{n}\left( L_{v_{j}}+L_{v_{j+1}}\right)  & \text{if }n%
\text{ is even} \\ 
&  \\ 
0_{D_{G}} & \text{if }n\text{ is odd.}
\end{array}
\right. $

\strut 

$\ $(5.18)$\ \ \ \ \ \ \ \ =\left\{ 
\begin{array}{ll}
2\mu _{n}\cdot 1_{D_{G}} & \text{if }n\text{ is even} \\ 
&  \\ 
0_{D_{G}} & \text{if }n\text{ is odd.}
\end{array}
\right. $

\strut 

Unfortunately, it is very hard to compute $\mu _{n},$ when $n\rightarrow
\infty .$ But we have to remark that if we have arbitrary graph $H$ and its
graph $W^{*}$-probability space $\left( W^{*}(H),F\right) $ over its
diagonal subalgebra $D_{H}$ and if $w\in loop^{c}(G),$ then

\strut 

$\ \ \ \ \ \ \ \ \ \mu _{w,w,...,w,w}^{1,*,...,1,*}=\mu _{n}=\mu
_{w,w,...,w,w}^{*,1,...,*,1},$ for all $n\in 2\mathbb{N}.$

\strut 

Now, let's compute the trivial $n$-th $D_{G}$-valued moment of $T.$ Notice
that since all odd trivial $D_{G}$-valued cumulants of $T$ vanish, all odd
trivial $D_{G}$-valued moments of $T$ vanish (See [11] and [14]). Thus it
suffices to compute the even trivial $D_{G}$-valued moments of $T.$ Assume
that $n\in 2\mathbb{N}.$ Then

\strut 

(5.19)\ \ \ \ \ \ \ \ \ $\ E\left( T^{n}\right) =\underset{\pi \in NC_{E}(n)%
}{\sum }k_{\pi }\left( T,...,T\right) ,$

\strut 

where $k_{\pi }\left( T,...,T\right) $ is the partition-dependent cumulant
of $T$ (See [16]) and

\strut 

$\ \ \ \ NC_{E}(n)\overset{def}{=}\{\pi \in NC(n):V\in \pi \Leftrightarrow
\left| V\right| \in 2\mathbb{N}\}.$

\strut 

By (5.18), we can get that $k_{n}(T,...,T)$ commutes with all elements in $%
W^{*}(G),$ because $1_{D_{G}}$ and $0_{D_{G}}$ commutes with $W^{*}(G)$ and $%
2\mu _{n}\in \mathbb{C},$ for all $n\in \mathbb{N}.$ So, the formula (5.19)
can be reformed by

\strut 

$\ \ \ E(T^{n})=\underset{\pi \in NC_{E}(n)}{\sum }\,\left( \underset{V\in
\pi }{\prod }k_{\left| V\right| }(T,...,T)\right) $

\strut 

$\ \ \ \ \ \ \ \ \ \ \ \ =\underset{\pi \in NC_{E}(n)}{\sum }\,\left( 
\underset{V\in \pi }{\prod }2\mu _{\left| V\right| }\cdot 1_{D_{G}}\right) $

\strut 

(5.20)$\ \ \ \ \ =\left( \underset{\pi \in NC_{E}(n)}{\sum }\,\left( 
\underset{V\in \pi }{\prod }2\mu _{\left| V\right| }\right) \right) \cdot
1_{D_{G}},$

\strut 

for all $n\in 2\mathbb{N}.$ Therefore, by (5.18) and (5.20), we have that if 
$T$ is the generating operator of the graph $W^{*}$-algebra of the circulant
graph $G$ with $N$-vertices, then

\strut \strut \strut 

$\ \ \ E(T^{n})=\left\{ 
\begin{array}{ll}
\left( \underset{\pi \in NC_{E}(n)}{\sum }\,\left( \underset{V\in \pi }{%
\prod }2\mu _{\left| V\right| }\right) \right) \cdot 1_{D_{G}} & \text{if }n%
\text{ is even} \\ 
&  \\ 
0_{D_{G}} & \text{if }n\text{ is odd.}
\end{array}
\right. $

and

$\ \ \ \ \ \ \ \ \ k_{n}\left( \underset{n\text{-times}}{\underbrace{%
T,.....,T}}\right) =\left\{ 
\begin{array}{ll}
2\mu _{n}\cdot 1_{D_{G}} & \text{if }n\text{ is even} \\ 
&  \\ 
0_{D_{G}} & \text{if }n\text{ is odd.}
\end{array}
\right. $

\strut \strut $\ $
\end{example}

\strut \strut \strut

\strut \strut

\begin{quote}
\textbf{Reference}

\strut

\strut

{\small [1] \ \ A. Nica, R-transform in Free Probability, IHP course note,
available at www.math.uwaterloo.ca/\symbol{126}anica.}

{\small [2]\strut \ \ \ A. Nica and R. Speicher, R-diagonal Pair-A Common
Approach to Haar Unitaries and Circular Elements, (1995), www
.mast.queensu.ca/\symbol{126}speicher.\strut }

{\small [3] \ }$\ ${\small B. Solel, You can see the arrows in a Quiver
Operator Algebras, (2000), preprint}

{\small \strut [4] \ \ A. Nica, D. Shlyakhtenko and R. Speicher, R-cyclic
Families of Matrices in Free Probability, J. of Funct Anal, 188 (2002),
227-271.}

{\small [5] \ \ D. Shlyakhtenko, Some Applications of Freeness with
Amalgamation, J. Reine Angew. Math, 500 (1998), 191-212.\strut }

{\small [6] \ \ D.Voiculescu, K. Dykemma and A. Nica, Free Random Variables,
CRM Monograph Series Vol 1 (1992).\strut }

{\small [7] \ \ D. Voiculescu, Operations on Certain Non-commuting
Operator-Valued Random Variables, Ast\'{e}risque, 232 (1995), 243-275.\strut 
}

{\small [10]\ D. Shlyakhtenko, A-Valued Semicircular Systems, J. of Funct
Anal, 166 (1999), 1-47.\strut }

{\small [10]\ D.W. Kribs and M.T. Jury, Ideal Structure in Free Semigroupoid
Algebras from Directed Graphs, preprint}

{\small [10]\ D.W. Kribs and S.C. Power, Free Semigroupoid Algebras, preprint%
}

{\small [11]\ I. Cho, Amalgamated Boxed Convolution and Amalgamated
R-transform Theory, (2002), preprint.}

{\small [12] I. Cho, The Tower of Amalgamated Noncommutative Probability
Spaces, (2002), Preprint.}

{\small [13] I. Cho, Free Perturbed R-transform Theory, (2003), Preprint.}

{\small [14]\ I. Cho, Compatibility of a Noncommutative Probability Space
and a Noncommutative Probability Space with Amalgamation, (2003), Preprint. }

{\small [15] I. Cho, Graph }$W^{*}${\small -Probability Spaces Over the
Diagonal Subalgebras, (2004), Preprint.}

{\small [16] R. Speicher, Combinatorial Theory of the Free Product with
Amalgamation and Operator-Valued Free Probability Theory, AMS Mem, Vol 132 ,
Num 627 , (1998).}

{\small [17] R. Speicher, Combinatorics of Free Probability Theory IHP
course note, available at www.mast.queensu.ca/\symbol{126}speicher.\strut }

{\small [18] T. Bates and D. Pask, Flow Equivalence of Graph Algebras,
(2004), Preprint.}
\end{quote}

\end{document}